\documentclass[12pt]{article}
\usepackage{amsmath,amssymb}
\usepackage[dvips]{graphicx,color}
\usepackage{a4wide}
\usepackage[authoryear]{natbib}
\usepackage{mathrsfs}
\usepackage{bm}
\usepackage{mathrsfs}

\numberwithin{equation}{section}
\renewcommand{\baselinestretch}{1.2}
\begin{document}
%\input D:/myfile/latex/Comlst
%\input Comlst
%%E:/zhengshurong2011/FmatrixApplications/S-1/Comlst
\newcommand{\ba}{{\bf a}}
\newcommand{\bA}{{\bf A}}
\newcommand{\bb}{{\bf b}}
\newcommand{\bB}{{\bf B}}
\newcommand{\bc}{{\bf c}}
\newcommand{\bC}{{\bf C}}
\newcommand{\bd}{{\bf d}}
\newcommand{\bD}{{\bf D}}
\newcommand{\bF}{{\bf F}}
\newcommand{\bx}{{\bf x}}
\newcommand{\bX}{{\bf X}}
\newcommand{\by}{{\bf y}}
\newcommand{\bY}{{\bf Y}}
\newcommand{\bz}{{\bf z}}
\newcommand{\bZ}{{\bf Z}}
\newcommand{\bw}{{\bf w}}
\newcommand{\bW}{{\bf W}}
\newcommand{\bI}{{\bf I}}
\newcommand{\bS}{{\bf S}}
\newcommand{\bT}{{\bf T}}
\newcommand{\balpha}{\bm{\alpha}}
\newcommand{\eprf}{\hspace*{\fill}$\blacksquare$}
\newcommand{\tr}{{\bf tr}}
\newcommand{\cov}{\mathop{\text cov}}

\newtheorem{thm}{Theorem}[section]
\newtheorem{prop}{Proposition}[section]
\newtheorem{lem}{Lemma}[section]
\newtheorem{rem}{Remark}[section]
\newtheorem{cor}{Corollary}[section]
\newtheorem{exm}{Example}[section]

%###################
\newcounter{theorem}
\newtheorem{theorem}{Theorem}[section]
\newcounter{lemma}
\newtheorem{lemma}{Lemma}[section]
\newcounter{remark}
\newtheorem{remark}{Remark}[section]
\newcounter{example}
\newtheorem{example}{Example}[section]
\newcounter{definition}
\newtheorem{definition}{Definition}[section]
\newcounter{corollary}
\newtheorem{corollary}{Corollary}[section]
\newcounter{proposition}
\newtheorem{proposition}{Proposition}[section]

\newcommand{\be}{\begin{equation}}
\newcommand{\ee}{\end{equation}}
\newcommand{\bqa}{\begin{eqnarray}}
\newcommand{\eqa}{\end{eqnarray}}
\newcommand{\bqn}{\begin{eqnarray*}}
\newcommand{\eqn}{\end{eqnarray*}}
\newcommand{\bdes}{\begin{description}}
\newcommand{\edes}{\end{description}}
\newcommand{\bitem}{\begin{itemize}}
\newcommand{\eitem}{\end{itemize}}
\newcommand{\bnum}{\begin{enumerate}}
\newcommand{\enum}{\end{enumerate}}
\newcommand{\ep}{\varepsilon}
\newcommand{\tep}{\tilde\varepsilon}
\newcommand{\hep}{\hat\varepsilon}
\newcommand{\ga}{\alpha}
\newcommand{\gb}{\beta}
\newcommand{\noin}{\noindent}
\newcommand{\DD}{\Delta_n}
\newcommand{\gl}{\lambda}
\newcommand{\gL}{\Lambda}
\newcommand{\Gl}{\Lambda}
\newcommand{\gd}{\delta}
\newcommand{\tgl}{\tilde\lambda}
\newcommand{\gs}{\sigma}
\newcommand{\gD}{\Delta}
\newcommand{\gG}{\Gamma}
\newcommand{\gS}{\Sigma}
\newcommand{\go}{\omega}
\newcommand{\gO}{\Omega}
\newcommand{\bXi}{{\mbox{\boldmath$\Xi$}}}
\newcommand{\bbxi}{{\mbox{\boldmath$\xi$}}}
\newcommand{\bbphi}{{\mbox{\boldmath$\phi$}}}
\newcommand{\bbPhi}{{\mbox{\boldmath$\Phi$}}}
\newcommand{\bbep}{{\mbox{\boldmath$\ep$}}}
\newcommand{\bgS}{{\mbox{\boldmath$\Sigma$}}}
\newcommand{\bga}{{\mbox{\boldmath$\alpha$}}}
\newcommand{\bbgd}{{\mbox{\boldmath$\delta$}}}
\newcommand{\bbtau}{{\mbox{\boldmath$\tau$}}}
\newcommand{\bbPsi}{{\mbox{\boldmath$\Psi$}}}
\newcommand{\bgma}{{\mbox{\boldmath$\gamma$}}}
\newcommand{\bbgG}{{\mbox{\boldmath$\Gamma$}}}
\newcommand{\btau}{{\mbox{\boldmath$\tau$}}}
\newcommand{\brho}{{\mbox{\boldmath$\rho$}}}
\newcommand{\bbgs}{{\mbox{\boldmath$\gs$}}}
\newcommand{\bbeta}{{\mbox{\boldmath$\beta$}}}
\newcommand{\bbzeta}{{\mbox{\boldmath$\zeta$}}}
\newcommand{\bbpsi}{{\mbox{\boldmath$\psi$}}}
\newcommand{\bbtheta}{{\mbox{\boldmath$\theta$}}}
\newcommand{\bdeta}{{\mbox{\boldmath$\eta$}}}
\newcommand{\bgl}{{\mbox{\boldmath$\lambda$}}}
\newcommand{\bbell}{{\mbox{\boldmath$\ell$}}}
\newcommand{\tell}{\tilde\ell}
\newcommand{\tbbell}{\widetilde{\mbox{\boldmath$\ell$}}}
\newcommand{\bzeta}{{\mbox{\boldmath$\zeta$}}}
\newcommand{\rtr}{{\rm tr}}
\newcommand{\Cov}{{\rm Cov}}
\newcommand{\rk}{{\rm rank}}
\newcommand{\sgn}{{\rm sgn}}
\newcommand{\ctan}{{\rm ctan}}
\newcommand{\rank}{{\rm rank}}
\newcommand{\rE}{{\rm E}}
\newcommand{\RRe}{{\rm Re}}
\newcommand{\IIm}{{\rm Im}}
\newcommand{\txi}{\tilde\xi}
\newcommand{\bxi}{{\mbox{\boldmath$\xi$}}}
\newcommand{\fii}{\frac{1}{2}}
\newcommand{\fnn}{\frac1{\sqrt n}}
\newcommand{\fn}{\frac{1}{n}}
\newcommand{\siln}{\sum_{i=1}^n}
\newcommand{\skln}{\sum_{k=1}^n}
\newcommand{\nskln}{\frac{1}{n}\sum_{k=1}^n}
\newcommand{\nsjln}{\frac{1}{n}\sum_{j=1}^n}
\newcommand{\sjln}{\sum_{j=1}^n}
\newcommand{\sklp}{\sum_{k=1}^p}
\newcommand{\psklp}{\frac{1}{p}\sum_{k=1}^p}
\newcommand{\nin}{\not\in}
\newcommand{\non}{\nonumber\\}
\newcommand{\lemm}{{\sc Lemma\ }}
\newcommand{\rema}{{\sc Remark\ }}
\newcommand{\proof}{{\sc Proof.\ \ }}
\newcommand{\ox}{{\overline x}}
\newcommand{\oy}{{\overline y}}
\newcommand{\hx}{{\widehat x}}
\newcommand{\oF}{{\overline F}}
\newcommand{\oX}{{\overline X}}
\newcommand{\obX}{{\overline \bbX}}
\newcommand{\oW}{{\overline \bXi}}
\newcommand{\omu}{{\overline \mu}}
\newcommand{\bbmu}{{\mbox{\boldmath$\mu$}}}
\newcommand{\bbnu}{{\mbox{\boldmath$\nu$}}}
\newcommand{\bpsi}{{\mbox{\boldmath$\Psi$}}}
\newcommand{\bsigma}{{\mbox{\boldmath$\sigma$}}}
\newcommand{\oox}{{\overline{\overline x}}}
\newcommand{\ooW}{{\overline{\overline \bXi}}}
\newcommand{\oomu}{{\overline{\overline \mu}}}
\newcommand{\tomu}{{\widetilde{\overline \mu}}}
\newcommand{\tox}{{\widetilde{\overline x}}}
\newcommand{\toXi}{{\widetilde{\overline \bXi}}}
\newcommand{\darrow}{\stackrel{\cal D}\to}
\def\iparrow{\buildrel i.p.\over\longrightarrow}
\def\asarrow{\buildrel a.s.\over\longrightarrow}
\newcommand{\tb}{{\tilde \bbb}}
\newcommand{\tx}{{\tilde x}}
\newcommand{\tnab}{{\tilde \nabla}}
\newcommand{\tXi}{{\tilde \bXi}}
\newcommand{\diag}{{\rm diag}}
\newcommand{\rP}{{\rm P}}
\newcommand{\rVar}{{\rm Var}}
\newcommand{\rCov}{{\rm Cov}}
\newcommand{\bbA}{{\bf A}}
\newcommand{\bba}{{\bf a}}
\newcommand{\bbB}{{\bf B}}
\newcommand{\bbb}{{\bf b}}
\newcommand{\bbC}{{\bf C}}
\newcommand{\bbc}{{\bf c}}
\newcommand{\bbD}{{\bf D}}
\newcommand{\bbd}{{\bf d}}
\newcommand{\bbe}{{\bf e}}
\newcommand{\bbE}{{\bf E}}
\newcommand{\bbf}{{\bf f}}
\newcommand{\bbF}{{\bf F}}
\newcommand{\bbg}{{\bf g}}
\newcommand{\bbG}{{\bf G}}
\newcommand{\bbH}{{\bf H}}
\newcommand{\bbh}{{\bf h}}
\newcommand{\bbI}{{\bf I}}
\newcommand{\bbi}{{\bf i}}
\newcommand{\bbj}{{\bf j}}
\newcommand{\bbJ}{{\bf J}}
\newcommand{\bbk}{{\bf k}}
\newcommand{\bbl}{{\bf 1}}
\newcommand{\bbL}{{\bf L}}
\newcommand{\bbM}{{\bf M}}
\newcommand{\bbm}{{\bf m}}
\newcommand{\bbn}{{\bf n}}
\newcommand{\bbN}{{\bf N}}
\newcommand{\bbQ}{{\bf Q}}
\newcommand{\bbp}{{\bf p}}
\newcommand{\bbP}{{\bf P}}
\newcommand{\bbq}{{\bf q}}
\newcommand{\bbO}{{\bf O}}
\newcommand{\bbR}{{\bf R}}
\newcommand{\bbr}{{\bf r}}
\newcommand{\bbs}{{\bf s}}
\newcommand{\bbS}{{\bf S}}
\newcommand{\tbbS}{\widetilde{\bf S}}
\newcommand{\hbbS}{\widehat{\bf S}}
\newcommand{\obbS}{\overline{\bf S}}
\newcommand{\bbT}{{\bf T}}
\newcommand{\bbt}{{\bf t}}
\newcommand{\hbbT}{\widehat{\bf T}}
\newcommand{\tbbT}{\widetilde{\bf T}}
\newcommand{\obT}{{\overline{\bf T}}}
\newcommand{\bbU}{{\bf U}}
\newcommand{\bbu}{{\bf u}}
\newcommand{\bbV}{{\bf V}}
\newcommand{\bbv}{{\bf v}}
\newcommand{\tbbv}{\widetilde{\bf v}}
\newcommand{\bbw}{{\bf w}}
\newcommand{\bbW}{{\bf W}}
\newcommand{\tbbW}{\widetilde{\bf W}}
\newcommand{\hbbB}{\widehat{\bf B}}
\newcommand{\hbbW}{\widehat{\bf W}}
\newcommand{\bbX}{{\bf X}}
\newcommand{\tbbB}{\widetilde {\bf B}}
\newcommand{\tbbX}{\widetilde {\bf X}}
\newcommand{\hbbX}{\widehat {\bf X}}
\newcommand{\bbx}{{\bf x}}
\newcommand{\obbx}{{\overline{\bf x}}}
\newcommand{\tbbx}{{\widetilde{\bf x}}}
\newcommand{\hbbx}{{\widehat{\bf x}}}
\newcommand{\obby}{{\overline{\bf y}}}
\newcommand{\hbbY}{{\widehat{\bf Y}}}
\newcommand{\tbbY}{{\widetilde{\bf Y}}}
\newcommand{\tbby}{{\tilde{\bf y}}}
\newcommand{\bbY}{{\bf Y}}
\newcommand{\bby}{{\bf y}}
\newcommand{\bbZ}{{\bf Z}}
\newcommand{\bbz}{{\bf z}}
\newcommand{\cC}{{\cal C}}
\newcommand{\cD}{{\cal D}}
\newcommand{\cE}{{\cal E}}
\newcommand{\cF}{{\cal F}}
\newcommand{\cG}{{\cal G}}
\newcommand{\cH}{{\cal H}}
\newcommand{\cI}{{\cal I}}
\newcommand{\cM}{{\cal M}}
\newcommand{\cP}{{\cal P}}
\newcommand{\cQ}{{\cal Q}}
\newcommand{\cR}{{\cal R}}
\newcommand{\cS}{{\cal S}}
\newcommand{\cT}{{\cal T}}
\newcommand{\cX}{{\cal X}}
\newcommand{\cU}{{\cal U}}
\newcommand{\cV}{{\cal V}}
\newcommand{\Th}{{\bbT^{1/2}}}
\newcommand{\um}{{\underline m}}
\newcommand{\us}{{\underline s}}
\newcommand{\umn}{{\underline m}_n}
\newcommand{\ua}{{\underline a}}
\newcommand{\uF}{{\underline F}}
\newcommand{\gma}{\gamma}
\newcommand{\bay}{\begin{array}}
\newcommand{\eay}{\end{array}}
\newcommand{\cas}[1]{\left\{\bay{ll}#1\eay\right.}
\newcommand{\matrx}[4]{\left[\begin{array}{cc}#1&#2\\ #3&#4\end{array}\right]}

\title{\vskip -2cm \bf CLT for  linear spectral statistics of
  random matrix $\bS^{-1}\bT$}
\author{\small
Shurong Zheng$^*$, Zhidong Bai$^*$ and Jianfeng Yao$^{**}$\\
\small $^*$School of Mathematics and Statistics and KLAS,
\small Northeast Normal University,\\
\small Changchun City 130024, P. R. China\\
\small $^{**}$Department of Statistics and Actuarial Science, Hong Kong
University, P. R. China}
\date{\small \today}
%\begin{document}
\renewcommand{\baselinestretch} {1.0}
\maketitle
\renewcommand{\baselinestretch} {1.5}
\begin{abstract}
As a generalization of the univariate Fisher statistic, random
Fisher matrices are widely-used  in multivariate statistical
analysis, e.g. for testing the equality of two multivariate
population covariance matrices. The  asymptotic distributions of
several meaningful test statistics depend on the related Fisher
matrices. Such Fisher matrices have the form ${\bf F}={\bf
S}_y{\bf M}{\bf S}_x^{-1}{\bf M}^*$ where $M$ is a non-negative
and non-random Hermitian matrix, and  ${\bf S}_x$ and ${\bf S}_y$
are $p\times p$ sample covariance matrices from two independent
samples where the populations are assumed centred and normalized
(i.e. mean 0, variance 1 and with independent components). In the
large-dimensional context, \citet{Zheng12}  establishes a central
limit theorem  for  linear spectral statistics of a standard
Fisher matrix where the two population covariance matrices are
equal, i.e. the matrix ${\mathbf M}$ is  the identity matrix and
${\bf F}={\bf S}_y{\bf S}_x^{-1}$. It is however of significant
importance to obtain a CLT for general Fisher matrices ${\mathbf
F}$ with an arbitrary ${\mathbf M}$ matrix. For the mentioned test
of equality, null distributions of test statistics rely on a
standard Fisher matrix with ${\mathbf M}=I_p$ while under the
alternative hypothesis, these distributions  depends on a general
Fisher matrix with arbitrary ${\mathbf M}$. As a first step to
this goal, we propose in this paper a CLT for spectral statistics
of the random matrix ${\bf S}_x^{-1}{\bf T}$ for a general
non-negative definite and {\bf non-random} Hermitian matrix $\bT$
(note that $\bT$ plays the role of ${\bf M}^*{\bf M}$). When ${\bf
T}$ is inversible, such a CLT can be directly derived using the
CLT of \citet{BS04} for the matrix ${\bT}^{-1}{\bf S}_x$. However,
in many large-dimensional statistic problems, the deterministic
matrix ${\bf T}$ is usually not inversible or has eigenvalues
close to zero.  The CLT from this paper covers this general
situation.
\end{abstract}

\section{Introduction}

For  a $p\times p$ random matrix $A_n$ with eigenvalues
$(\lambda_j)$,  linear spectral statistics (LSS)
of type  $\frac1p \sum_j f(\lambda_j)$ for various test functions $f$
are of  central importance in the theory of random matrices
and its applications
Central limit theorems (CLT) for such LSS
of large dimensional random matrices
have a long history,
and received considerable attention in recent years.
They have
important applications in various domains like number theory,
high-dimensional multivariate statistics and wireless communication
networks; for more information, the readers are referred
to the recent survey paper \citet{John07}.
To mention a few, in an early work,
\citet{Jons82}
 gave a  CLT for $(\rtr(\bA_n), \cdots,\rtr(\bA_n^k))$ for a sequence of
 Wishart matrices $(\bA_n)$, where $k$ is a fixed number, and the dimension $p$ of the matrices grows
proportionally to the sample size $n$. Subsequent works include
\citet{CL95}, \citet{JOH} which considered  extensions of
classical Gaussian ensembles, and \citet{SS98a,SS98b}  where
Gaussian fluctuations are identified for LSS of Wigner matrices
with a class of more general test  functions. A general CLT for
LSS of Wigner matrices was given  in \citet{BY05} where in
partiular, the limiting mean and covariance functions are
identified. Similarly,  \citet{BS04} established a CLT  for
general sample covariance matrices with explicit limiting
parameters. In \citet{LytPas09}, the authors reconsider such CLTs
but with a new idea of interpolation  that allows the
generalisation from Gaussian matrix ensembles to matrix ensembles
with general entries satisfying a moment condtiion. Recent
improvments are proposed in \citet{PanZhou08} that propose a
generalization of the CLT in \cite{BS04} (see also \citet{WY13}
for a complement on these CLT's). Finally, \citet{Pan12} and
\citet{BZ13} extend \cite{BS04}'s CLT to biased and unbiased
sample covariance matrices, respectively.

Random Fisher matrices are widely-used  in multivariate
statistical analysis, e.g. for testing the equality of two
multivariate population covariance matrices. The  asymptotic
distributions of several meaningful test statistics depend on the
related Fisher matrices. Such Fisher matrices have the form ${\bf
F}={\bf S}_y{\bf M}{\bf S}_x^{-1}{\bf M}^*$ where ${\bf M}$ is a
non-negative deterministic Hermitian matrix, and  ${\bf S}_x$ and
${\bf S}_y$ are $p\times p$ sample covariance matrices from two
independent  samples where the populations are assumed centred and
normalized (i.e. mean 0, variance 1 and with independent
components). In the large-dimesional context, \citet{Zheng12}
establishes a CLT   for  linear spectral statistics of a standard
Fisher matrix where the two population covariance matrices are
equal, i.e. the matrix ${\mathbf M}$ is  the identity matrix and
${\bf F}={\bf S}_y {\bf S}_x^{-1}$. It is however of significant
importance to obtain a CLT for general Fisher matrices ${\mathbf
F}$ with an arbitrary ${\mathbf M}$ matrix. For the mentioned test
of equality,  null distributions of test statistics rely  on a
standard Fisher matrix with ${\mathbf M}=I_p$ while under the
alternative hypothesis, these distributions  depends on a general
Fisher matrix with arbitrary ${\mathbf M}$.

In order to extend the CLT of \citet{Zheng12} to general Fisher
matrices,  we first need to establish limit theorems for the
spectral (eigenvalues)  distribution of the matrix ${\bf M}{\bf
S}_x^{-1}{\bf M}^*$, or the matrix ${\bf S}_x^{-1}{\bf T}$ where
$\bT ={\bf M}^*{\bf M}$ is {\bf non-random}. This includes i) an
identifiation of the limit of its  spectral distribution; ii) a
CLT for its LSS. When the {\bf non-random} matrix ${\bf T}$ is
inversible, since ${\bf S}_x^{-1}{\bf T}= \left[
\bT^{-1}\bS_x\right]^{-1}$, CLT for LSS of  ${\bf S}_x^{-1}{\bf
T}$ can be derived from the CLT of \citet{BS04}. However, in many
large-dimensional statistic problems, the deterministic  matrix
$\bT$ is usually not invertible or has eigenvalues close to zero,
and it is then hopeless to base the analysis on the  CLT of
\cite{BS04}.

In this paper, we consider the product  ${\bf S}_x^{-1}{\bf T}$ of
a general determinist and {\bf non-random} Hermitian matrix $\bbT$
by the inverse ${\bf S}_x^{-1}$ of a standard sample covariance
matrix. As the main results of the paper, solutions to the
aforementioned problems are provided.

The organization of this paper is as follows.
Section~\ref{sec:main} presents our main results.
The proofs of these two main theorems are given in the following
sections,
respectively.

\section{Main results}
\label{sec:main}

Following \cite{BS04}, let
$\{\bx_t\}$, $t=1,\ldots,n$ be  a sequence of indepenent
$p$-dimensional observations
with independent and standardised components, i.e. for
$\bx_t=(x_{tj})$,
$\rE x_{tj}=0$ and $\rE|x_{tj}|^2=1$. The corresponding sample covariance
matrix is
\begin{equation}
  \label{eq:S}
  \bS=\frac1n\sum_{t=1}^n \bx_t \bx_t^*~.
\end{equation}
Consider the product matrix
\begin{equation}
  {\bf S}^{-1}{\bf T}=\left(\frac{1}{n}\sum\limits_{t=1}^n
  {\bf x}_t{\bf x}^*_t\right)^{-1}{\bf T}~,
\end{equation}
where $\bT$ is a $p\times p$ non-negative definite and non-random
Hermitian matrix. Notice that we do not ask  $\bT$ be invertible.

We first state the  framework for our main results.
\begin{description}
\item[{\bf Assumption 1}]
  The $p\times n$ observation matrix
  $(x_{tj}, t=1,\cdots,n, j=1,\cdots,p)$
  are made with independent elements satisfying
  $\rE x_{tj}=0$, $\rE
 |x_{tj}|^2=1$.  Moreover,
 for any $\eta>0$ and as $p,~n\to\infty$,
 \begin{equation}
   \frac{1}{np}\sum\limits_{t=1}^n\sum\limits_{j=1}^p
   \rE\left[|x_{tj}|^2I_{\{|x_{tj}|\geq\eta\sqrt{n}\}}\right]\rightarrow
   0~,
   \label{eq:Lindeberg2}
 \end{equation}
 where $I_{\{\cdot\}}$ is the indicator function.

 The elements are either all real or all complex and
 we  set an index  $\kappa=1$ or $\kappa=2$, respectively. In the
 later case,   $\rE  \{x_{tj}^2\}=0$ for all $t,j$.

\item[{\bf Assumption 1$^{*}$}]
  In addition to Assumption 1, the entries $\{x_{tj}\}$ have an
  uniform 4-th moment
  $\rE |x_{tj}|^4=1+\kappa$.  Moreover,
  for any $\eta>0$  and as $p,~n\to\infty$,
  \begin{equation}
    \frac{1}{np}\sum\limits_{t=1}^n\sum\limits_{j=1}^p
    \rE \left[|x_{tj}|^4I_{\{|x_{tj}|\geq\eta\sqrt{n}\}}\right]\rightarrow
0.
   \label{eq:Lindeberg4}
  \end{equation}

\item[{\bf Assumption 1$^{**}$}]
  In addition to Assumption 1, the entries $\{x_{tj}\}$ have
  a finite
  4-th moment  (not necessarily the same).
  Moreover,
  for any $\eta>0$  and as $p,~n\to\infty$,
  \begin{equation}
    \frac{1}{np}\sum\limits_{t=1}^n\sum\limits_{j=1}^p
    \rE\left[|x_{tj}|^4I_{\{|x_{tj}|\geq\eta\sqrt{n}\}}\right]\rightarrow
0.
   \label{eq:Lindeberg4}
  \end{equation}

\item [{\bf Assumption 2}]
  The ESD $H_n$ of  $\{\bbT\}$ tends to a limit
  $H$, which is a  probability measure not degenerated to
  the Dirac mass at 0.
  %\footnote{previously    said:      ``un-degenerated at 0''. is it OK?}.
\item[{\bf Assumption 2$^*$}]
  In addition to Assumption 2,
  the operator norm of $\bbT$ is bounded when $n,~p\to\infty$.
\item[{\bf Assumption 3}]
  The dimension $p$ and the sample size $n$ both tend to infinity
  such that $p/n\to y\in (0,1)$.
\end{description}

Assumption 1 states that the entries are independent, not necessarily
identically distributed, but with homogeneous moments of first and
second order, together with a Lindeberg type condition of order
2.
 Assumption 1$^*$ reinforce Assumption 1 with similar condtions
using a homogeneous  forth order moment that matches
the Gaussian case.
 Assumption 1$^{**}$ generalizes the
previous one by allowing arbitrary values for the fourth moment
of the entries.

Recall that the empirical spectral distribution (ESD) of a matrix is
the distribution generated by its eigenvalues.
When this ESD has a limit when the dimensions grow to infinity, the
limit
is called the limiting spectral distribution (LSD) of the matrix.

The first main result of the paper identifies the LSD
of $\bS^{-1}\bT$.

\begin{thm}\label{thm1}
  Under Assumptions 1, 2 and 3, with probability 1, the ESD  $F_n$ of
  $\bbS^{-1}\bbT$ tends to a non-random distribution $F^{y,H}$ whose
  Stieltjes
  transform $s(z)$ is   the unique solution to the equation
{\color{red}\be zs(z)=-1+\int\frac{tdH(t)}{-z-yz^2s(z)+t}
~.\label{eqsilv3} \ee} The distribution  $F^{y,H}$ is then  the
LSD of $\bbS^{-1}\bbT$.
\end{thm}

Next, we consider a   LSS of  $\bbS^{-1}\bbT$ of form
\[
 F_n (f) = \int f(x)dF_n (x) = \frac1p \sum_{j=1}^p f(\lambda_j)~,
\]
where the $\{\lambda_j\}$'s are the eigenvalues of the matrix
$\bbS^{-1}\bbT$ and $f$ a given test function.
Similarly to \cite{BS04}, a special feature here is that
fluctuations of $F_n(f)$ will  not be considered around
the LSD limit $F^{y,H}(f)$, but around
$F^{y_n,H_n}(f)$,
a finite-sample  proxy of $F^{y,H}$ obtained
by substituting the parameters
$(y_n,H_n)$ to
$(y,H)$  in the LSD. Therefore, we consider the random variable
\[
X_n(f)=p\left[  F_n(f) -  F^{y_n,H_n}(f) \right]
= p \int f(x) d[F_n-F^{y_n,H_n}](x)~.
\]

The second main result of the paper is the following CLT.

\begin{thm}\label{thm2} Assume that
  Assumptions 1$^*$, 2$^*$ and 3 hold.
  Let $f_1,\cdots,f_k$ be functions analytic on an open domain
  of the complex plane enclosing the interval
  \[
  \left[\frac{\lim\inf\limits_{p}\lambda_{\min}^\bT}{(1+\sqrt{y})^2},
    ~~\frac{\lim\sup\limits_{p}\lambda_{\max}^\bT}{(1-\sqrt{y})^2}\right]~,
  \]
  where $\lambda_{\min}^\bT$ and $\lambda_{\max}^\bT$ are respectively the
  smallest and the largest eigenvalue of $\bT$.
  Then,
  the random vector
  $ [X_n( f_1),\cdots,X_n( f_k)]$
  weakly  converges to a Gaussian vector
  $[X_{f_1},\cdots,X_{f_k}]$ with mean function
  \begin{equation}\label{mean1}
    \begin{array}{lll}
      \rE X_{f_j}&=&\displaystyle{-\frac{\kappa-1}{2\pi i}\oint
        f_j(z)\cdot\frac{1}{z^2} \frac{y\int
          \frac{t(1+yzs(z))^3dH(t)}{(t/z-1-yzs(z))^3}}{\left(
          1-y\int\frac{(1+yzs(z))^2 dH(t)}{(t/z-1-yzs(z))^2}\right)^2}dz}~,
    \end{array}
  \end{equation}
  and
  covariance function
  \begin{equation}\label{cov1}
    Cov(X_{f_i},X_{f_j})=-\frac{\kappa}{4\pi^2}
    \oint\oint\frac{f_i(z_1)f_j(z_2)\cdot\frac{\partial[z_1(1+yz_1s(z_1))]}{\partial z_1}\frac{\partial[z_2(1+yz_2s(z_2))]}{\partial z_2}}{[z_1(1+yz_1s(z_1))-z_2(1+yz_2s(z_2))]^2}
    dz_1dz_2~
  \end{equation}
  where $\frac{1}{z^2}\cdot \frac{y\int
\frac{t(1+yzs(z))^3dH(t)}{(t/z-1-yzs(z))^3}}{\left(
1-y\int\frac{(1+yzs(z))^2
dH(t)}{(t/z-1-yzs(z))^2}\right)^2}=\frac12\frac{d\log\left(1-y\int\frac{(1+yzs(z))^2
dH(t)}{(t/z-1-yzs(z))^2}\right)}{dz}.$
  The contours in (\ref{mean1}) and (\ref{cov1}) are closed and are
  taken in the positive direction in the complex plane,
  all enclosing the support of $F^{y,H}$.
\end{thm}

When the fourth moments of the entries are different from the
value $\kappa+1$  matching the Gaussian case (3 or 2), the
expression (\ref{var}) has an additional term
\[
\frac{1}{n^2}\sum\limits_{i=1}^nb_i(z_1)b_i(z_2)
\sum\limits_{j=1}^p    (\rE|X_{ij}|^4-1-\kappa)
\left[\rE_{i-1}\left(z_1^{-1}\bbT-\bS_i\right)^{-1}\right]_{jj}
\left[\rE_{i-1}\left({z_2}^{-1}\bbT-\bS_i\right)^{-1}\right]_{jj}
\]
and the expression (\ref{mean}) has an additional term
\[
-\frac{\rE(\beta_1(z))^2}{n^2y}\sum\limits_{j=1}^p
(\rE|X_{1j}|^4-1-\kappa)
\left[\left(z^{-1}\bbT-\bS_1\right)^{-1}\right]_{jj}
\left[\left(z^{-1}\bbT-\bS_1\right)^{-1}\bbT\left(z^{-1}\bbT-E\beta_1(z)\bI\right)^{-1}\right]_{jj}.
\]
Then the covariance (\ref{var_complement}) and mean (\ref{bias})
will have additional terms, the limits of
$$
\frac{\partial^2\left\{\frac{1}{n^2}\sum\limits_{i=1}^nb_i(z_1)b_i(z_2)
\sum\limits_{j=1}^p    (\rE|X_{ij}|^4-1-\kappa)
\left[\rE_{i-1}\left(z_1^{-1}\bbT-\bS_i\right)^{-1}\right]_{jj}
\left[\rE_{i-1}\left({z_2}^{-1}\bbT-\bS_i\right)^{-1}\right]_{jj}\right\}}{\partial
z_1\partial z_2}
$$
and
$$
\begin{array}{l}
\frac{y(1+yzs(z))^3}{z^2\left( 1-y\int\frac{(1+yzs(z))^2
dH(t)}{(t/z-1-yzs(z))^2}\right)}\cdot\frac{1}{p}
\sum\limits_{j=1}^p\bigg\{(E|X_{1j}|^4-1-\kappa)\left[\left(z^{-1}\bbT-\bS_1\right)^{-1}\right]_{jj}
\quad\quad\quad\quad\quad\quad\\
\quad\quad\quad\quad\quad\quad\quad\quad\quad\quad\quad\quad\quad\quad\quad\quad\quad\quad\quad\quad\cdot\left[\left(z^{-1}\bbT-\bS_1\right)^{-1}\bbT\left(z^{-1}\bbT-E\beta_1(z)\bI\right)^{-1}\right]_{jj}\bigg\}.
\end{array}
$$
When $\frac{1}{p}\sum\limits_{j=1}^p(E|X_{ij}|^4-1-\kappa)
    \left[\left({z_1}^{-1}\bbT-\bS_i\right)^{-1}\right]_{jj}
    \left[\left({z_2}^{-1}\bbT-\bS_i\right)^{-1}\right]_{jj}$
    converges to $h(z_1,z_2)$ uniformly in $i$, then the covariance
    (\ref{var_complement}) will have the additional term
$$
\frac{\partial^2[y\cdot(1+yz_1s(z_1))(1+yz_2s(z_2))h(z_1,z_2)]}{\partial
z_1\partial z_2}
$$
because $\rE b_i(z)\rightarrow 1+yzs(z)$ by (\ref{7}).

Then Theorem~\ref{thm2} is easily  extended to this situation as
follows.

\begin{prop}\label{prop}
Assume that
  Assumptions 1$^{**}$, 2$^{*}$ and 3 hold.
  Let $f_1,\cdots,f_k$ be functions analytic on an open domain
  of the complex plane enclosing the interval
  \[
  \left[\frac{\lim\inf\limits_{p}\lambda_{\min}^\bT}{(1+\sqrt{y})^2},
    ~~\frac{\lim\sup\limits_{p}\lambda_{\max}^\bT}{(1-\sqrt{y})^2}\right]~,
  \]
  where $\lambda_{\min}^\bT$ and $\lambda_{\max}^\bT$ are respectively the
  smallest and the largest eigenvalue of $\bT$. Moreover,
  assume in addition that the following  non-random limits exist:
  \begin{enumerate}
  \item
    $\frac{1}{p}\sum\limits_{j=1}^p(E|X_{ij}|^4-1-\kappa)
    \left[\left({z_1}^{-1}\bbT-\bS_i\right)^{-1}\right]_{jj}
    \left[\left({z_2}^{-1}\bbT-\bS_i\right)^{-1}\right]_{jj}$
    converges to $h(z_1,z_2)$ uniformly in $i$;
  \item
    $\frac{1}{p}\sum\limits_{j=1}^p(E|X_{1j}|^4-1-\kappa)\left[\left(z^{-1}\bbT-\bS_1\right)^{-1}\right]_{jj}
    \left[\left(z^{-1}\bbT-\bS_1\right)^{-1}\bbT\left(z^{-1}\bbT-E\beta_1(z)\bI\right)^{-1}\right]_{jj}$
    converges to $h_M(z)$.
  \end{enumerate}
  Then the random vector
  $ [X_n( f_1),\cdots,X_n( f_k)]$
  weakly  converges to a Gaussian vector
  $[X_{f_1},\cdots,X_{f_k}]$ with mean function
 \begin{eqnarray}
 \rE X_{f_j}&=&-\frac{\kappa-1}{2\pi i}\oint
        f_j(z)\cdot\frac{1}{z^2} \frac{y\int
          \frac{t(1+yzs(z))^3dH(t)}{(t/z-1-yzs(z))^3}}{\left(
          1-y\int\frac{(1+yzs(z))^2
          dH(t)}{(t/z-1-yzs(z))^2}\right)^2}dz,\nonumber\\
           &&-\frac{1}{2\pi i}\int
      f_j(z)\cdot\frac{y(1+yzs(z))^3}{z^2\left(1-y\int\frac{(1+yzs(z))^2dH(t)}{(t/z-1-yzs(z))^2}\right)}\cdot
      h_M(z)dz\nonumber
 \end{eqnarray}
 and covariance function
 \begin{eqnarray}
 Cov(X_{f_i},X_{f_j})&=&-\frac{\kappa}{4\pi^2}
    \oint\oint\frac{f_i(z_1)f_j(z_2)\cdot
    \frac{\partial[z_1(1+yz_1s(z_1))]}{\partial z_1}\frac{\partial[z_2(1+yz_2s(z_2))]}
    {\partial z_2}}{[z_1(1+yz_1s(z_1))-z_2(1+yz_2s(z_2))]^2}
    dz_1dz_2\nonumber\\
    &&-\frac{1}{4\pi^2}\oint\oint f_i(z_1)f_j(z_2)\cdot
    \frac{\partial^2[y\cdot(1+yz_1s(z_1))(1+yz_2s(z_2))h(z_1,z_2)]}{\partial z_1\partial
    z_2}dz_1dz_2.\nonumber
     \end{eqnarray}
\end{prop}
The contours are closed and are
  taken in the positive direction in the complex plane,
  all enclosing the support of $F^{y,H}$.

\vskip 0.5cm When $E|X_{ij}|^4-1-\kappa=\beta_x+o(1)$ uniformly in
$i, j$ and $\bbT$ is a diagonal matrix with positive eigenvalues,
then we have
\begin{eqnarray}
&&\frac{1}{p}\sum\limits_{j=1}^p(E|X_{1j}|^4-1-\kappa)
\left[\left(z^{-1}\bbT-\bS_1\right)^{-1}\right]_{jj}
\left[\left(z^{-1}\bbT-\bS_1\right)^{-1}\bbT\left(z^{-1}\bbT-E\beta_1(z)\bI\right)^{-1}\right]_{jj}\nonumber\\
&\rightarrow&h_M(z)=\beta_x\cdot\int\frac{z^3t
dH(t)}{[t-z(1+yzs(z))]^3}\nonumber
\end{eqnarray}
and
\begin{eqnarray}
&&\frac{1}{p}\sum\limits_{j=1}^p(\rE|X_{ij}|^4-1-\kappa)
\left[\left(z_1^{-1}\bbT-\bS_i\right)^{-1}\right]_{jj}
\left[\left({z_2}^{-1}\bbT-\bS_i\right)^{-1}\right]_{jj}\nonumber\\
&=&h(z_1,z_2)=\beta_x\cdot \int\frac{z_1z_2
dH(t)}{[t-z_1(1+yz_1s(z_1))][t-z_2(1+yz_2s(z_2))]}.\nonumber
\end{eqnarray}
Then the mean (\ref{bias}) has the additional term
$$
\frac{\beta_x\cdot yz(1+yzs(z))^3}{\left(
1-y\int\frac{(1+yzs(z))^2 dH(t)}{(t/z-1-yzs(z))^2}\right)}\cdot
\int\frac{tdH(t)}{[t-z(1+yzs(z))]^3}
$$
and the covariance (\ref{var_complement}) has the additional term
$$
\beta_x\cdot\frac{\partial^2}{\partial z_1\partial
z_2}\left[y\int\frac{z_1(1+yz_1s(z_1))z_2(1+yz_2s(z_2))
}{[t-z_1(1+yz_1s(z_1))][t-z_2(1+yz_2s(z_2))]}dH(t)\right].
$$
Then Proposition \ref{prop} easily extends to the following
proposition.

\begin{prop}
Let assumptions of Proposition \ref{prop} hold. Moreover, assume
that $E|X_{ij}|^4-1-\kappa=\beta_x+o(1)$ uniformly in $i,j$ and
$\bbT$ is a diagonal matrix with positive eigenvalues, then we
obtain that $ [X_n( f_1),\cdots,X_n( f_k)]$
  weakly  converges to a Gaussian vector
  $[X_{f_1},\cdots,X_{f_k}]$ with mean function
 \begin{eqnarray}
 \rE X_{f_j}&=&-\frac{\kappa-1}{2\pi i}\oint
        f_j(z)\cdot\frac{1}{z^2} \frac{y\int
          \frac{t(1+yzs(z))^3dH(t)}{(t/z-1-yzs(z))^3}}{\left(
          1-y\int\frac{(1+yzs(z))^2
          dH(t)}{(t/z-1-yzs(z))^2}\right)^2}dz,\nonumber\\
           &&-\frac{\beta_x}{2\pi i}\int
     \left[ f_j(z)\cdot\frac{yz(1+yzs(z))^3}{1-y\int\frac{(1+yzs(z))^2 dH(t)}{(t/z-1-yzs(z))^2}}
      \int\frac{tdH(t)}{[t-z(1+yzs(z))]^3}\right]dz\nonumber
 \end{eqnarray}
 and covariance function
 \begin{eqnarray}
 &&Cov(X_{f_i},X_{f_j})\nonumber\\
 &=&-\frac{\kappa}{4\pi^2}
    \oint\oint\frac{f_i(z_1)f_j(z_2)\cdot
    \frac{\partial[z_1(1+yz_1s(z_1))]}{\partial z_1}\frac{\partial[z_2(1+yz_2s(z_2))]}
    {\partial z_2}}{[z_1(1+yz_1s(z_1))-z_2(1+yz_2s(z_2))]^2}
    dz_1dz_2\nonumber\\
    &&-\frac{\beta_x}{4\pi^2}\oint\oint\left\{f_i(z_1)f_j(z_2)\cdot
\frac{\partial^2}{\partial z_1\partial
z_2}\left[y\int\frac{z_1(1+yz_1s(z_1))z_2(1+yz_2s(z_2))
dH(t)}{[t-z_1(1+yz_1s(z_1))][t-z_2(1+yz_2s(z_2))]}\right]\right\}
    dz_1dz_2.\nonumber
     \end{eqnarray}
The contours in (\ref{mean1}) and (\ref{cov1}) are closed and are
  taken in the positive direction in the complex plane,
  all enclosing the support of $F^{y,H}$.

\end{prop}

\subsection{Relation between Theorem \ref{thm2}  and
  the CLT in  \protect\citet{BS04}}

Theorem~\ref{thm2} can be viewed a complemet to the CLT in
\citet{BS04} while moving from
the sample covariance matrix $\bS$ to its
inverse $\bS^{-1}$.
When the factor $\bT$ in $\bS^{-1}\bT$ is not invertible, these
CLT's are not directly comparable.
If $\bT$ is indeed invertible,
these CLT's should  be comparable.
In this subsection, we will prove that they are indeed the same in
this case.
More precisely we prove that the
mean and covariance functions given in Theorem~\ref{thm2}
are the same as those given in
Theorem 1.1 of \citet{BS04}.

Actually,  when ${\bf T}$
is invertible, we have \bqn
s_n(z)&=&\frac1p\rtr(\bS^{-1}\bT-z\bI)^{-1}=\frac1p\sum\limits_{i=1}^p\frac{1}{\lambda_i((\bS\bT^{-1})^{-1})-z}\\
      &=&\frac1p\sum\limits_{i=1}^p\frac{1}{1/\lambda_i(\bS\bT^{-1})-z}
      =\frac1p\sum\limits_{i=1}^p\frac{\lambda_i(\bS\bT^{-1})}{1-z\lambda_i(\bS\bT^{-1})}\\
      &=&\frac{-1}{pz}\sum\limits_{i=1}^p\frac{\lambda_i(\bS\bT^{-1})}{\lambda_i(\bS\bT^{-1})-\frac{1}{z}}\\
      &=&-\frac{1}{z}-\frac{1}{z^2}\cdot\frac{1}{p}
      \sum\limits_{i=1}^p\frac{1}{\lambda_i(\bS\bT^{-1})-\frac{1}{z}}\\
      &=&-\frac{1}{z}-\frac{1}{z^2}\cdot s_n^{\bS\bT^{-1}}(\frac1z)\\
      &=&-\frac{1}{z}-\frac{1}{z^2}\cdot\bigg(\frac1y\underline{m}_n(\frac1z)+\frac{1-y}{y}z\bigg)\\
      &=&-\frac{1}{yz}-\frac{1}{yz^2}\underline{m}_n(\frac1z)~,
\eqn
where $\underline{m}_n$ is the Stieltjes transform of
$\bX_n\bT^{-1}\bX_n^{*}$ with $\bX_n^*=(\bbx_1,\cdots,\bbx_n)$ is
$p\times n$. That is,
\begin{equation}\label{sz}
s(z)=-\frac{1}{yz}-\frac{1}{yz^2}\cdot\underline{m}(\frac1z),\quad
1+yzs(z)=1-1-\frac{1}{z}\underline{m}(\frac1z)=-\frac{1}{z}\underline{m}(\frac1z)~,
\end{equation}
where $s(z)$ is the limit of $s_n(z)$ and $\underline{m}(z)$ is
the limit of $\underline{m}_n(z)$. So the CLT of $p(s_n(z)-s(z))$
is the same as $-\frac{1}{z^2}\cdot
n(\underline{m}_n(\frac1z)-\underline{m}_n(\frac1z))$. By Lemma
1.1 of Bai and Silverstein (2004), we know that the CLT of
$$
-\frac{1}{z^2}\cdot n(\underline{m}_n(\frac1z)-\underline{m}_n(\frac1z))
~,$$
has mean
$$
-\frac{1}{z^2}\frac{y\int\frac{t\cdot (\underline{m}(1/z))^3d(H(t))}{(t+\underline{m}(1/z))^3}}{[1-y\int\frac{(\underline{m}(1/z))^2d(H(t))}
{(t+\underline{m}(1/z))^2}]^2}~,
$$
and covariance
$$
\frac{1}{z_1^2z_2^2}\frac{\underline{{m}}'(\frac{1}{z_1})\underline{{m}}'(\frac{1}{z_2})}
{(\underline{{m}}(\frac{1}{z_2})-\underline{{m}}(\frac{1}{z_1}))^2}
-\frac{1}{(z_1-z_2)^2}~.
$$
It is easily to verify that
$$
-\frac{1}{z^2}\frac{y\int\frac{t\cdot (\underline{m}(1/z))^3d(H(t))}{(t+\underline{m}(1/z))^3}}{[1-y\int\frac{(\underline{m}(1/z))^2d(H(t))}
{(t+\underline{m}(1/z))^2}]^2}
=\frac{1}{z^2}\frac{y\int\frac{t\cdot z^3(1+ys(z))^3dH(t)}{(t-z(1+yzs(z)))^3}}{[1-y\int\frac{z^2(1+yzs(z))^2dH(t)}{(t-z(1+yzs(z)))^2}]^2}
=\frac{1}{z^2}\frac{y\int\frac{t\cdot(1+ys(z))^3dH(t)}{(t/z-1-yzs(z))^3}}{[1-y\int\frac{(1+yzs(z))^2dH(t)}{(t/z-1-yzs(z))^2}]^2}~,
$$
and
$$
\frac{1}{z_1^2z_2^2}\frac{\underline{{m}}'(\frac{1}{z_1})\underline{{m}}'(\frac{1}{z_2})}
{(\underline{{m}}(\frac{1}{z_2})-\underline{{m}}(\frac{1}{z_1}))^2}
-\frac{1}{(z_1-z_2)^2}=\frac{[z_1(1+yz_1s(z_1))]'[z_2(1+yz_2s(z_2))]'}{[z_1(1+yz_1s(z_1))-z_2(1+yz_2s(z_2))]^2}
-\frac{1}{(z_1-z_2)^2}~,
$$
which are the same as given in Theorem~\ref{thm2}. Thus, when ${\bf T}$ is
inversible, the CLT of LSS of $\bbS^{-1}\bbT$ has the same mean
and covariance functions as that obtained by Theorem 1.1 of Bai
and Silverstein (2004).
%\section{Comments and Conclusions}

%\section*{Acknowledgement}
%The authors wish to thank an editor, an associate editor, and two
%referees for their valuable comments and constructive suggestions
%which greatly improve the presentation of the paper. Zheng's
%research was partially supported by NSFC-11171058 and
%NECT-11-0616.

%
\section{Proof of Theorem \ref{thm1}}

Using exactly the same approach employed in Section 4.3 of Bai and
Silverstein (2010), we may truncate the extreme eigenvalues of
$\bbT$ and tails of the random  variables $x_{ij}$ and then
renormalize them without altering the LSD of $\bbS^{-1}\bbT$. So
we may assume that Assumption 2$^*$ is true and
$|x_{ij}|\le\eta_n\sqrt{n}$ where $\eta_n\to 0$.

Now, we proceed with the proof of Theorem \ref{thm1}. To start
with, we assume that $\bbT$ is invertible and there is a positive
constant $\omega>0$ such that $H(\omega)=0$, that is, the norm of
$\bbT^{-1}$ is bounded. By Theorem 4.1 of Bai and Silverstein
(2010) we know that the LSD of $\bbS\bbT^{-1}$ exists and its
Stieltjes transform $m(z)$ satisfies \be
m(z)=\int\frac{1}{t(1-y-yzm(z))-z}dH(1/t)=\int\frac{tdH(t)}{1-y-yzm(z)-tz}.
\label{eqsilv} \ee Note that $m(z)$ is the unique solution to the
equation (\ref{eqsilv}) that has the same sign of  imaginary part
as $z$.

If we denote the Stieltjes transforms of $\bbS^{-1}\bbT$ and $\bbS\bbT^{-1}$ by $s_n(z)$ and $m_n(z)$, respectively. By the relation
$$
m_n(z)=-\frac1z-\frac{1}{z^2}s_n(1/z),
$$
and $m_n(z)\to m(z)$ a.s., we know that with probability 1,
$s_n(z)$ converges to a limit $s(z)$ that satisfies \be
-\frac1z-\frac{1}{z^2}s(1/z)=\int\frac{tdH(t)}{1-y-yz(-\frac1z-\frac{1}{z^2}s(1/z))-tz}.
\label{eqsilv2} \ee Changing $z$ as $1/z$ and simplifying it, we
obtain (\ref{eqsilv3}).

Now, we consider possibly singular $\bT$ and  will
show that for any fixed $z=u+iv$ with $v>0$, $s_n(z)$ still  converges
to a limit $s(z)$ that satisfies (\ref{eqsilv3}).

For any fixed $\ep>0$, define $\bbT_\ep=\bbT+\ep\bbI$ and define $\bbS_+$ from $\bbS$ by replacing its eigenvalues less than $\frac12a$ as $\frac12 a$,
where $a=(1-\sqrt{y})^2$. By the rank inequality of Bai (1999), we have
\be
\left\|F^{\bbS^{-1}\bbT}-F^{\bbS^{-1}_+\bbT}\right\|\le \frac1p\ ^{\#}\left\{\lambda_i(\bbS)\le \frac12a\right\}\to 0, a.s.
\label{eqb1}
\ee
By Theorem A.45 of Bai and Silverstein (2010),
\be
L(F^{\bbS^{-1}_+\bbT},F^{\bbS^{-1}_+\bbT_\ep})\le \|\bbS^{-1}_+(\bbT-\bbT_\ep)\|\le 2a^{-1}\ep.
\label{eqb2}
\ee
Using again the rank inequality, we have
\be
\left\|F^{\bbS^{-1}\bbT_\ep}-F^{\bbS^{-1}_+\bbT_\ep}\right\|\le \frac1p\ ^{\#}\left\{\lambda_i(\bbS)\le \frac12a\right\}\to 0, a.s.
\label{eqb3}
\ee
By what has been proved anove for invertible $\bT$,
with probability 1, $s_{n,\ep}(z)=\frac1p\rtr(\bbS^{-1}\bbT_\ep)\to s_\ep(z)$ which is a solution to the equation
\be
zs_\ep(z)=-1+\int\frac{tdH_\ep(t)}{1-z-yz^2s_\ep(z)+t}.
\label{eqsilv4}
\ee
where $H_\ep(t)=H(t-\ep)$.

To complete the proof of the theorem, we only need to verify that
the equation (\ref{eqsilv4}) has a unique solution that is the
Stieltjes transform of a probability distribution, and the
solution $s_\ep(z)$ is right-continuous at $\ep=0$. Making a
transformation $w_\ep(z)=\sqrt{z}(1+z s_\ep(z)$, where $\sqrt{z}$
is the square root of $z$ satisfying $\Im(z)\Im(\sqrt{z})>0$, then
the equation (\ref{eqsilv4}) becomes \be
w_\ep(z)=\int\frac{tdH_\ep(t)}{\frac{1+t}{\sqrt{z}}-(1-y)\sqrt{z}-yw_\ep(z)},
\label{eqsilv5} \ee where $w_\ep(z)$ has the same sign of imaginary
part as $z$.

We only need to consider the case where $\Im(z)>0$. Let
$w_2=\Im(w_\ep(z))>0$, comparing the imaginary parts of
(\ref{eqsilv5}), we have \bqn w_2&=&\int
\frac{\frac{(1+t)\Im(\sqrt{z})}{|z|}+(1-y)\Im(\sqrt{x})+yw_2}
{\left|\frac{1+t}{\sqrt{z}}-(1-y)\sqrt{z}-yw_\ep(z)\right|^2}dH_\ep(t)\\
&>&\int \frac{yw_2}
{\left|\frac{1+t}{\sqrt{z}}-(1-y)\sqrt{z}-yw_\ep(z)\right|^2}dH_\ep(t),
\eqn which implies that \be \int \frac{ydH_\ep(t)}
{\left|\frac{1+t}{\sqrt{z}}-(1-y)\sqrt{z}-yw_\ep(z)\right|^2}<1.
\label{eqbb1} \ee

Suppose (\ref{eqsilv5}) had two solution $w^{(j)}$
with $w^{(j)}_2=\Im(w^{(j)})>0$, $j=1,2$. Then making difference of both
sides and cancelling $w_1-w_2$ from both sides, we obtain
$$
1=y\int\frac{tdH_\ep(t)}{(\frac{1+t}{\sqrt{z}-(1-y)\sqrt{z}-yw^{(1)}})
(\frac{1+t}{\sqrt{z}-(1-y)\sqrt{z}-yw^{(2)}})}~,
$$
which implies by Cauchy-Schwarz that
$$
1\le \left(\int \frac{ydH_\ep(t)}
{\left|\frac{1+t}{\sqrt{z}}-(1-y)\sqrt{z}-yw^{(1)}\right|^2}\int \frac{ydH_\ep(t)}
{\left|\frac{1+t}{\sqrt{z}}-(1-y)\sqrt{z}-yw^{(2)}\right|^2}\right)^{1/2}<1,
$$
where the last inequality follows by applying (\ref{eqbb1}) for both
$w^{(1)}$ and $w^{(2)}$. The contradiction proves the uniqueness of a
solution to (\ref{eqsilv5}).

Finally, we show that the solution $w_\ep$ is right-continuous at $\ep=0$.  By (\ref{eqsilv5}), we have
\bqa
w_\ep(z)-w_0(z)&=&\frac{\int\frac{td(H_\ep(t)-H(t))}
{\frac{1+t}{\sqrt{z}}-(1-y)\sqrt{z}-yw_\ep(z)}}{1-y\int\frac{tdH(t)}
{\left(\frac{1+t}{\sqrt{z}}-(1-y)\sqrt{z}-yw_\ep(z)\right)
\left(\frac{1+t}{\sqrt{z}}-(1-y)\sqrt{z}-yw_0(z)\right)}}.
\label{eqbb2}
\eqa
Since
\bqn
&&\left|\frac{1+t}{\sqrt{z}}-(1-y)\sqrt{z}-yw_\ep(z)\right|\ge
-\Im(\frac{1+t}{\sqrt{z}}-(1-y)\sqrt{z}-yw_\ep(z))>(1-y)\Im(\sqrt{z})
~,
\eqn
we have
\bqn
&&\left|y\int\frac{tdH(t)}
{\left(\frac{1+t}{\sqrt{z}}-(1-y)\sqrt{z}-yw_\ep(z)\right)
\left(\frac{1+t}{\sqrt{z}}-(1-y)\sqrt{z}-yw_0(z)\right)}\right|\\
&\le& \left(y\int\frac{tdH(t)}
{\left|\frac{1+t}{\sqrt{z}}-(1-y)\sqrt{z}-yw_0(z)\right|^2}\right)^{1/2}<1.
\eqn
It follows that $w_\ep(z)-w_0(z)\to 0$ which implies that
$s_\ep(z)-s(z)\to 0$.

 The proof of Theorem \ref{thm1} is complete.

\section{Proof of Theorem \ref{thm2}}

We first describe the
strategy of the proof that follows
the proof in \citet{BS04}
and an improved version in
\citet{BS10}.
First, due to  Assumption 1$^*$,
we may truncate the random variables $x_{ij}$ at $\eta_n\sqrt{n}$
and renormalize them without alerting the CLT of $X_n(f)$, where
$\eta_n\downarrow 0$ with some slow rate. Therefore, we may make
the following additional assumptions: \bnum \item $|x_{ij}|\le
\eta_n\sqrt{n}$; \item $\rE x_{ij}^2=\kappa-1+o(n^{-1})$; \item
$\rE |x_{ij}^4|=1+\kappa+o(1)$. \enum Define a contour $\cC_n$ by
$$
{\cal{C}}_n={\cal{C}}_l\cup{\cal{C}}_u\cup{\cal{C}}_b\cup{\cal{C}}_r
$$
where
$$
 \begin{array}{ll}
 {\cal{C}}_u=\{x+i\nu_0: x\in[x_l,x_u]\},&
 {\cal{C}}_b=\{x-i\nu_0: x\in[x_l,x_u]\},\\
 {\cal{C}}_l=\{x_l+i\nu: |\nu|\le \nu_0\},&
 {\cal{C}}_r=\{x_r+i\nu: |\nu|\le \nu_0\}
\end{array}
$$
and $(x_l,x_r) \supset [\liminf
\gl_{\min}(\bbT)/(1+\sqrt{y})^2,\limsup
\gl_{\max}(\bbT)/(1-\sqrt{y})^2]$ and is enclosed in the analytic
region of the $f_j(x)$'s. Following Bai and Silverstein (2004), we
can rewrite $X_n(f)$ as \be X_b(f)=-\frac1{2\pi}\oint_{\cC_n}f(z)
p(s_n(z)-s_n^0(z)) dz, \label{remk} \ee where $s_n^0(z)$ is the
Stieltjes transform of $F^{y_n,H_n}$.
\begin{remark}
Note that the identity (\ref{remk}) holds only when all
eigenvalues of $\bbS^{-1}\bbT$ are falling inside the interval
$(x_l,x_r)$. By Bai and Yin (1993), with probability 1, when $n$
is large, all eigenvalues of $\bbS$ are falling inside the
interval $((1-\sqrt{y})^2-\ep,(1+\sqrt{y})^2+\ep)$ which confirms
that (\ref{remk}) holds for all large $n$. So, without loss of
generality, in the proof of Theorem \ref{thm2}, we assume
(\ref{remk}) holds.
\end{remark}

Write $M_n(z)=p(s_n(z)-s_n^0(z))=M_n^{1}(z)+ M_n^{2}(z)$, where
$M_n^{1}(z)= p(s_n(z)-\rE s_{n}(z))$ and $M_n^{2}(z)=p(\rE
s_n(z)-s_n^0(z))$. We shall establish a CLT for $M_n^{1}(z)$,  and
 then find the limit of $M_n^{2}(z)$ on $\cC_u$ and $\cC_b$.
Their combinaison will  complete
the proof of Theorem \ref{thm2}.

\subsection{Finite-dimensional convergence of $M_n^{1}(z)$ on $\cC_u$}
%Because the support set of LSD of $\bS^{-1}\bT$ is as follows
%$$
%\left[\frac{\lim\inf\limits_{p}\lambda_{min}^\bT}{(1+\sqrt{y})^2}, \frac{\lim\sup\limits_{p}\lambda_{max}^\bT}{(1-\sqrt{y})^2}\right],
%$$
%where $\lambda_{min}^\bT$ and $\lambda_{max}^\bT$ are the minimum and maximum eigenvalues of
%$\bT$, then let $\nu_0>0$ be arbitrary, $x_l$ be a real number smaller than $\frac{\lim\inf\limits_{p}\lambda_{min}^\bT}{(1+\sqrt{y})^2}$
%and $x_r$ be a real number greater than $\frac{\lim\sup\limits_{p}\lambda_{max}^\bT}{(1-\sqrt{y})^2}$. Then

%Let $M_n(z)=M_n^1(z)+M_n^2(z)$ where
%$$
%M_n(z)=p(s_n(z)-s(z)),\quad M_n^1(z)=p(s_n(z)-\rE s_n(z))\quad\mbox{and}\quad M_n^2(z)=p(\rE s_n(z)-s(z)).
%$$
%Let the truncated version of $M_n(z)$ be
%$$
%\hat{M}_n(z)=\left\{
%\begin{array}{ll}
%M_n(z),&z\in{\cal{C}}_n\\
%M_n(x_r+in^{-1}\epsilon_n),& \nu\in[0,n^{-1}\epsilon_n]\\
%M_n(x_l+in^{-1}\epsilon_n),& \nu\in[0,n^{-1}\epsilon_n]
%\end{array}
%\right.
%$$
%and $\hat{M}_n^1(z)$ can be defined similarly.
We first prove an auxiliary theorem.
\begin{thm}\label{thm3}
Under Assumptions 1$^*$, 2$^*$, and 3, $M_n^{1}(z)$ converges
weakly to a complex Gaussian process $M_1(\cdot)$ on the contour
$z\in\mathcal{C}$, with mean function
$$
\rE M_1(z)=0
$$
and covariance function
\begin{equation}\label{eqn1}
\Cov(M_1(z_1),
M_1(z_2))=\kappa\left[\frac{[z_1(1+yz_1s(z_1))]'[z_2(1+yz_2s(z_2))]'}{[z_1(1+yz_1s(z_1))-z_2(1+yz_2s(z_2))]^2}
-\frac{1}{(z_1-z_2)^2}\right].
\end{equation}
\end{thm}
\proof Let $\rE_i$ denote the conditional expectation given
$\{\bbx_1,\cdots,\bbx_i\}$ and $\rE_0$ denote the unconditional
expectation. Denote $z=u+iv$ with $v>0$ fixed,
$$
\begin{array}{ll}
\balpha_i=\frac{1}{\sqrt{n}}{\bf X}_i,~{\bf S}_k={\bf
S}-\balpha_k\balpha_k^*,&\bD=\bbS-z\bbI,~\bD_k=\bT-z\bS_k,\\
\beta_i(z)=\frac{1}{1-\balpha_i'(\frac{1}{z}\bT-\bS_i)^{-1}\balpha_i},&
\bar\beta_i(z)=\frac{1}{1-\frac{1}{n}\rtr(\frac{1}{z}\bT-\bS_i)^{-1}},\\
\bar{\beta}_i(z,\theta_i)=\frac{1}{1-\frac{1}{n}\rtr(\frac{1}{z}\bT-\bS_i)^{-1}+\theta_i\hat{\gamma}_i(z)},&
\hat\gamma_i(z)=\frac{1}{n}\rtr(\frac{1}{z}\bT-\bS_i)^{-1}-\balpha_i'(\frac{1}{z}\bT-\bS_i)^{-1}\balpha_i.
\end{array}
$$
 Then we have
 $$
\beta_i(z)^{-1}=\bar{\beta}_i(z)^{-1}-\hat{\gamma}_i(z),\quad
\bar\beta_i(z)\beta_i^{-1}(z)=1-\bar{\beta}_i(z)\hat{\gamma}_i(z).
$$ Therefore, by Taylor expansion \be
\begin{array}{lll}
&&(\rE_{i}-\rE_{i-1})\log\beta_i^{-1}(z)\\
&=&(\rE_{i}-\rE_{i-1})\left(\log\beta^{-1}_i(z)-\log\bar\beta^{-1}_i(z)\right)\\
&=&(\rE_{i}-\rE_{i-1})[-\bar{\beta}_i(z)\hat{\gamma}_i(z)
+\bar{\beta}_i^2(z,\theta_i)\hat{\gamma}_i^2(z)]\\
&=&-\rE_{i}\bar{\beta}_i(z)\hat{\gamma}_i(z)+(\rE_{i}-\rE_{i-1})\bar{\beta}_i^2(z,\theta_i)\hat{\gamma}_i^2(z).
\end{array}
\label{eqzhuzhi}
\ee
Here, we have used a formula that $\log \beta_i^{-1}-\log \bar\beta_i^{-1}
=\log \bar\beta_i(z)\beta_i^{-1}$. In fact we should add an additional
term $2\pi k(z)$ where $k(z)$ is a random integer function of
$z$. This term does make any contribution because we only need the
derivative  of the function $\log \beta_i^{-1}$ in the next step.

For any $i\le n$, we have
$$
\begin{array}{lll}
s_n(z)&=&\frac1p\rtr(\bS^{-1}\bT-z\bI)^{-1}=\frac1p\rtr\bS(\bT-z\bS)^{-1}\\
      &=&\frac1p\rtr(\bS_i+\balpha_i\balpha_i')(\bT-z\bS_i-z\balpha_i\balpha_i')^{-1}\\
      &=&\frac{1}{pz}\rtr(z\bS_i+z\balpha_i\balpha_i')\left((\bT-z\bS_i)^{-1}+\frac{z(\bT-z\bS_i)^{-1}
      \balpha_i\balpha_i'(\bT-z\bS_i)^{-1}}{1-z\balpha'_i(\bT-z\bS_i)^{-1}\balpha_i}\right)\\
      &=&-\frac{1}{pz}\rtr(\bT-z\bS_i-z\balpha_i\balpha_i'-\bT)\left((\bT-z\bS_i)^{-1}+\frac{z(\bT-z\bS_i)^{-1}
      \balpha_i\balpha_i'(\bT-z\bS_i)^{-1}}{1-z\balpha'_i(\bT-z\bS_i)^{-1}\balpha_i}\right)\\
      &=&-\frac{1}{pz}\left(p+\frac{z\balpha_i'\bD_i^{-1}\balpha_i}{1-z\balpha'_i\bD_i^{-1}\balpha_i}
      -\frac{z^2(\balpha_i'\bD_i^{-1}\balpha_i)^2}{1-z\balpha'_i\bD_i^{-1}\balpha_i}
      -z\balpha_i'\bD_i^{-1}\balpha_i\right)\\
      &&+\frac{1}{pz}\rtr\left(\bT\bD_i^{-1}+\frac{z\bT\bD_i^{-1}
      \balpha_i\balpha_i'\bD_i^{-1}}{1-z\balpha'_i\bD_i^{-1}\balpha_i}\right)\\
      &=&-\frac{1}{z}+\frac{1}{pz}\rtr\bT\bD_i^{-1}+\frac{1}{p}\frac{
      \balpha_i'\bD_i^{-1}\bT\bD_i^{-1}\balpha_i}{1-z\balpha'_i\bD_i^{-1}\balpha_i}\\
      &=&-\frac{1}{z}+\frac{1}{pz}\rtr\bT\bD_i^{-1}
      +\frac1p\frac{\partial}{\partial z}\log\beta_i(z)~.
      \end{array}
      $$
Therefore,
$$
\begin{array}{lll}
s_n(z)-\rE s_n(z)&=&-\frac1p\frac{\partial}
{\partial z}\sum\limits_{i=1}^n(\rE_{i}-\rE_{i-1})\log(1-z\balpha'_i\bD_i^{-1}\balpha_i)\\
&=&-\frac1p\frac{\partial}
{\partial z}\sum\limits_{i=1}^n\rE_{i}
\bar{\beta}_i(z)\hat{\gamma}_i(z)-\frac1p\frac{\partial}
{\partial z}\sum\limits_{i=1}^n(\rE_i-\rE_{i-1})\bar{\beta}_i^2(z,\theta_i)\hat{\gamma}_i^2(z).
\end{array}
$$
Since
\begin{equation}
\begin{array}{lll}
&&E\left|\sum\limits_{i=1}^n(\rE_{i-1}-\rE_i)\bar{\beta}^2_i(z,\theta_i)\hat{\gamma}^2_i(z)\right|^2\\
&=&\sum\limits_{i=1}^nE\left|(\rE_{i-1}-\rE_i)\bar{\beta}^2_i(z,\theta_i)\hat{\gamma}^2_i(z)\right|^2\\
&\leq&4\sum\limits_{i=1}^nE\left|\bar{\beta}^2_i(z,\theta_i)\hat{\gamma}^2_i(z)\right|^2\\
&\leq&4\sum\limits_{i=1}^nE\left|\hat{\gamma}_i(z)\right|^4\\
&=&nO(n^{-1}\eta_n^4)=o(1)\ ~\mbox{(by Lemma \ref{BaiLem91})}
\end{array}~,\label{5}
\end{equation}
 we have
$$
p((s_n(z)-\rE s_n(z))=-\sum\limits_{i=1}^nE_i\frac{d}{dz}\bar\beta_i(z)\hat{\gamma}_i(z)
+o_p(1)=-\frac{d}{dz}\sum\limits_{i=1}^nY_i(z)+o_p(1)~,
$$
where
$$\begin{array}{lll}Y_i(z)&=&\rE_i\bar{\beta}_i(z)\hat{\gamma}_i(z).%\\
%&=&\rE_i\bar{\beta}_i^2(z)\frac1n\rtr(\bD_i^{-1}+z\bD_i^{-1}\bS_i\bD_i^{-1})\hat{\gamma}_i(z)\\
%&&+\bar{\beta}_i(z)
%\left(\frac1n\rtr(\bD_i^{-1}+z\bD_i^{-1}\bS_i\bD_i^{-1})-\balpha_i'(\bD_i^{-1}
%+z\bD_i^{-1}\bS_i\bD_i^{-1})\balpha_i\right)
\end{array}
$$
We first consider a finite  sum
$$
\sum\limits_{k=1}^r\sum\limits_{i=1}^na_kY_i(z_k)=\sum\limits_{i=1}^n\sum\limits_{k=1}^ra_kY_i(z_k)
$$
from  $r$ points $z_k$ on the contour with arbitrary weighting numbers $a_k$.
That is, we need to complete the following two steps:
\par
\noindent {\em Step 1:} Verify the Lyapunoff condition, i.e. $\sum\limits_{i=1}^n\rE
\left|Y_i(z)\right|^4=o(1).$

In fact, if $z\in \cC_u$ or $\cC_b$, by the fact that $|\bar\beta_i(z)|<|z|/\nu_0$,
\bqn
\sum\limits_{i=1}^n\rE
\left|Y_i(z)\right|^4&\le & C\sum\limits_{i=1}^n\rE
\left|\hat\gamma_i(z)\right|^4\\
&\le& Cn^{-4}\sum_{i=1}^n\rE\bigg[\left(\rtr(\bbD_i^{-1}(z)\bbD_i^{-1}(\bar z))\right)^2+\max_{ij}\rE|x_{ij}^8\sum_{j=1}^p|[\bbD_i^{-1}]_{jj}|^4\bigg]\\
&\le& C(n^{-1}+\eta_n^4)\to 0,
\eqn
where $[\bbD_i^{-1}]_{jj}$ is the $j$-th diagonal entry of $\bbD_i^{-1}(z)$ which is bounded by $|z|/\nu_0$.

\medskip\noindent
{\em Step 2:}  Find the limits of $\frac{\partial^2}{\partial z_1\partial
z_2}\sum\limits_{i=1}^n \rE
_{i-1}Y_i(z_1)Y_i(z_2)=\frac{\partial^2}{\partial z_1\partial
z_2}\sum\limits_{i=1}^n \rE _{i-1}
[\rE_i\bar{\beta}_i(z_1)\hat{\gamma}_i(z_1)\cdot
\rE_i\bar{\beta}_i(z_2)\hat{\gamma}_i(z_2)]$.

{\color{red}
$$
E\left|\sum\limits_{i=1}^n\sum\limits_{k=1}^ra_kY_i(z_k)\right|^2=
\sum\limits_{i=1}^nE\left|\sum\limits_{k=1}^ra_kY_i(z_k)\right|^2
\leq
K\sum\limits_{i=1}^n\sum\limits_{k=1}^r|a_k|^2E\left|Y_i(z_k)\right|^2\leq
K
$$
because $E\left|Y_i(z_k)\right|^2=O(n^{-1})$ by Lemma
\ref{BaiLem91}.} We have
$$
\begin{array}{lll}
\bar{\beta}_i(z)-b_i(z)&=&\bar{\beta}_i(z)b_i(z)\left(n^{-1}\rtr(\frac{1}{z}\bT-\bS_i)^{-1}-
n^{-1}\rE\rtr(\frac{1}{z}\bT-\bS_i)^{-1}\right)
\end{array}
$$
where $b_i(z)=\frac{1}{1-n^{-1}\rE\rtr(\frac{1}{z}\bT-\bS_i)^{-1}}$.
Then (Bai and Silverstein (2010), P139)
$$
\begin{array}{lll}
E|\bar{\beta}_i(z)-b_i(z)|^{2l}&\leq&
K\rE\left(n^{-1}\rtr(\frac{1}{z}\bT-\bS_i)^{-1}-
n^{-1}\rE\rtr(\frac{1}{z}\bT-\bS_i)^{-1}\right)^{2l}\\
&=&K\rE\left(\frac{1}{n}\sum\limits_{j=2}^{n}\left[\rE_j\rtr(\frac{1}{z}\bT-\bS_1)^{-1}-
\rE_{j-1}\rtr(\frac{1}{z}\bT-\bS_1)^{-1}\right]\right)^{2l}\\
&=&K\rE\Bigg(\frac{1}{n}\sum\limits_{j=2}^{n}
(\rE_j-\rE_{j-1})\rtr\Bigg\{(\frac{1}{z}\bT-\bS_1)^{-1}-(\frac{1}{z}\bT-\bS_{1j})^{-1}\Bigg\}
\Bigg)^{2l}\\
&=&K\rE\Bigg(\frac{1}{n}\sum\limits_{j=2}^{n}
(\rE_j-\rE_{j-1})\frac{\alpha_j'(\frac{1}{z}\bT-\bS_{1j})^{-2}\alpha_j}
{1+\alpha_j'(\frac{1}{z}\bT-\bS_{1j})^{-1}\alpha_j}
\Bigg)^{2l}\\
&\leq&\frac{K}{n^{2l}}\rE\Bigg(\sum\limits_{j=2}^{n}
\left|(\rE_j-\rE_{j-1})\frac{\alpha_j'(\frac{1}{z}\bT-\bS_{1j})^{-2}\alpha_j}
{1+\alpha_j'(\frac{1}{z}\bT-\bS_{1j})^{-1}\alpha_j}\right|^2\Bigg)^{l}\\
&&(\mbox{Lemma 2.12 of Bai and Silverstein (2010)})\\
&\leq&\displaystyle{\frac{K}{n^l\cdot\nu^{2l}}}.
\end{array}
$$
Then $\rE|\bar{\beta}_i(z)-b_i(z)|^{2l}=O(n^{-l})$ is uniformly.
Then
$$
\sum\limits_{i=1}^n\rE _{i-1}
[\rE_i\bar{\beta}_i(z_1)\hat{\gamma}_i(z_1)\cdot
\rE_i\bar{\beta}_i(z_2)\hat{\gamma}_i(z_2)]-\sum\limits_{i=1}^n
b_i(z_1)b_i(z_2)\rE _{i-1}[\rE_i\hat{\gamma}_i(z_1)\cdot
\rE_i\hat{\gamma}_i(z_2)]=o_p(1)
$$
Then we only consider the limit of
\begin{equation}\label{var}
\begin{array}{lll}
&&\sum\limits_{i=1}^nb_i(z_1)b_i(z_2)\rE_{i-1}[\rE_i\hat{\gamma}_i(z_1)\rE_i\hat{\gamma}_i(z_2)]\\
\\
&=&\frac{\kappa}{n^2}\sum\limits_{i=1}^nb_i(z_1)b_i(z_2)
\rtr\left[\rE_i(\frac{1}{z_1}\bT-\bS_i)^{-1}\rE_i(\frac{1}{z_2}\bT-\bS_i)^{-1}\right]\\
&&\mbox{[by (1.15) of Bai and Silverstein (2004)]}.
\end{array}
\end{equation}
We have
$$
(\frac{1}{z}\bT-\bS_i)-\frac{1}{z}\bT+\frac{n-1}{n}b_i(z)\bI=-\sum\limits_{k\not=i}\balpha_k\balpha_k^*
+\frac{n-1}{n}b_i(z)\bI.
$$
Multiplying by $(\frac{n-1}{n}b_i(z)\bI-\frac{1}{z}\bT)^{-1}$ on
the left, $(\frac{1}{z}\bT-\bS_i)^{-1}$ on the right, then we have
\bqa\label{inverse}
(\frac{1}{z}\bT-\bS_i)^{-1}&=&-(\frac{n-1}{n}b_i(z)\bI-\frac{1}{z}\bT)^{-1}\nonumber\\
                           &&-\sum\limits_{k\not=i}\beta_{k(i)}(z)
                           (\frac{n-1}{n}b_i(z)\bI-\frac{1}{z}\bT)^{-1}\balpha_k\balpha_k^*
                           (\frac{1}{z}\bT-\bS_{ik})^{-1}\nonumber\\
                           &&+\frac{n-1}{n}b_i(z)(\frac{n-1}{n}b_i(z)\bI-\frac{1}{z}\bT)^{-1}
                           (\frac{1}{z}\bT-\bS_i)^{-1}\nonumber\\
                           &=&-(\frac{n-1}{n}b_i(z)\bI-\frac{1}{z}\bT)^{-1}-b_i(z)\bA(z)
                           -\bB(z)-C(z)
\eqa where
$$
\bA(z)=\sum\limits_{k\not=i}
                           \bigg(\frac{n-1}{n}b_i(z)\bI-\frac{1}{z}\bT\bigg)^{-1}
                           (\balpha_k\balpha_k^*-\frac{1}{n}\bI)(\frac{1}{z}\bT-\bS_{ik})^{-1}
$$
$$
\bB(z)=\sum\limits_{k\not=i}
(\beta_{k(i)}(z)-b_i(z))\bigg(\frac{n-1}{n}b_j(z)\bI-\frac{1}{z}\bT\bigg)^{-1}
\balpha_k\balpha_k^*(\frac{1}{z}\bT-\bS_{ik})^{-1}
$$
$$
\bC(z)=\frac{1}{n}b_i(z)
\bigg(\frac{n-1}{n}b_i(z)\bI-\frac{1}{z}\bT\bigg)^{-1}
\sum\limits_{k\not=i}\bigg[(\frac{1}{z}\bT-\bS_{ik})^{-1}-(\frac{1}{z}\bT-\bS_{i})^{-1}\bigg]
$$
Similarly, we have
$$
\left\|\bigg(\frac{n-1}{n}b_j(z)\bI-\frac{1}{z}\bT\bigg)^{-1}\right\|\leq
K
$$
where $K$ is a constant. Similarly, we have \bqa\label{E2}
|b_{ik}(z)-b_i(z)|&=&\left|b_i(z)b_{ik}(z)
\left[\frac{1}{n}\rE\rtr(\frac{1}{z}\bT-\bS_i)^{-1}-\frac{1}{n}\rE\rtr(\frac{1}{z}\bT-\bS_{ik})^{-1}\right]\right|\nonumber\\
                  &=&\frac{1}{n}\left|b_{ik}(z)b_i(z)\rE\beta_{k(i)}(z)
\balpha_k^*(\frac{1}{z}\bT-\bS_{ik})^{-2}\balpha_k\right|=O(n^{-1})
\eqa and
\begin{equation}\label{E3}
\rE(\beta_{k(i)}(z)-b_{ik}(z))^2=O(n^{-1})~\mbox{by Lemma
\ref{BaiLem91}}
\end{equation}
where
$$
\beta_{k(i)}(z)=\frac{1}{1-\balpha_k'(\frac{1}{z}\bT-\bS_{ik})^{-1}\balpha_k}
$$
$$
b_i(z)=\frac{1}{1-\frac{1}{n}\rE\rtr(\frac{1}{z}\bT-\bS_i)^{-1}},~
b_{ik}(z)=\frac{1}{1-\frac{1}{n}\rE\rtr(\frac{1}{z}\bT-\bS_{ik})^{-1}}.
$$
First we have
$$
\begin{array}{lll}
&&\rE\left|\rtr\rE_i\bB(z_1)rE_i(\frac{1}{z_2}\bT-\bS_i)^{-1}\right|\\
&\leq&\sum\limits_{k\not=i}
\rE\left|\rE_i\rtr(\beta_{k(i)}(z_1)-b_i(z_1))(\frac{n-1}{n}b_j(z_1)\bI-\frac{1}{z_1}\bT)^{-1}
\balpha_k\balpha_k^*(\frac{1}{z_1}\bT-\bS_{ik})^{-1}(\frac{1}{z_2}\bT-\breve{\bS}_i)^{-1}\right|\\
&\leq&\sum\limits_{k\not=i}\rE\left|\rE_i(\beta_{k(i)}(z)-b_i(z))
\balpha_k^*(\frac{1}{z}\bT-\bS_{ik})^{-1}(\frac{1}{z_2}\bT-\breve{\bS}_i)^{-1}\balpha_k\right|\\
&\leq&\sum\limits_{k\not=i}
\rE^{1/2}|(\beta_{k(i)}(z)-b_i(z))^2|\cdot
\rE^{1/2}\left|\balpha_k^*(\frac{1}{z}\bT-\bS_{ik})^{-1}(\frac{1}{z_2}\bT-\breve{\bS}_i)^{-1}\balpha_k\right|\\
&=&\sqrt{n}\rE^{1/2}\left|\balpha_k^*(\frac{1}{z}\bT-\bS_{ik})^{-1}(\frac{1}{z_2}\bT-\breve{\bS}_i)^{-1}\balpha_k\right|
~\mbox{(by
(1.15) of Bai and
Silverstein (2004))}\\
&=&O(n^{1/2})
\end{array}
$$
where $\breve{\bS}_i$ is the analogue for the matrix $\bS_i$ with
vectors $\bx_{j+1},\cdots,\bx_n$ replaced by their iid copies
$\breve{\bx}_{j+1},\cdots,\breve{\bx}_n$.

Second we have
$$
\begin{array}{lll}
&&\frac1n\sum\limits_{k\not=i}\rE\left|\rtr\rE_i\bC(z_1)\cdot\rE_i
(\frac{1}{z_2}\bT-\bS_i)^{-1}\right|\\
&\leq&\frac{1}{n}\sum\limits_{k\not=i}\rE\left|\rtr\rE_ib_i(z)
(\frac{n-1}{n}b_j(z)\bI-\frac{1}{z}\bT)^{-1}
\bigg((\frac{1}{z_1}\bT-\bS_{ik})^{-1}-(\frac{1}{z}\bT-\bS_i)^{-1}\bigg)(\frac{1}{z_2}\bT-\breve{\bS}_i)^{-1}\right|\\
&=&\frac{1}{n}\sum\limits_{k\not=i}\rE\bigg|\rE_ib_i(z)
(\frac{n-1}{n}b_j(z)\bI-\frac{1}{z}\bT)^{-1}
\beta_{k(i)}\balpha_k'(\frac{1}{z_1}\bT-\bS_{ik})^{-1}(\frac{1}{z_2}\bT-\breve{\bS}_i)^{-1}
(\frac{1}{z_1}\bT-\bS_{ik})^{-1}\balpha_k\bigg|\\
&\leq&\frac{1}{n}\sum\limits_{k\not=i}
\rE\bigg|\balpha_k'(\frac{1}{z_1}\bT-\bS_{ik})^{-1}(\frac{1}{z_2}\bT-\breve{\bS}_i)^{-1}
(\frac{1}{z_1}\bT-\bS_{ik})^{-1}\balpha_k\bigg|\\
&&~\mbox{(by (1.15) of Bai and Silverstein (2004))}\\
&\leq& K
\end{array}
$$
Third, we consider
$$
\begin{array}{lll}
&&b_i(z_1)\rtr\rE_i\bA(z_1)\rE_i(\frac{1}{z_2}\bT-\bS_i)^{-1}\\
&=&b_i(z_1)\rtr
\sum\limits_{k<i}(\frac{n-1}{n}b_i(z_1)\bI-\frac{1}{z_1}\bT)^{-1}
(\balpha_k\balpha_k^*-\frac{1}{n}\bI)\rE_i(\frac{1}{z_1}\bT-\bS_{ik})^{-1}\rE_i(\frac{1}{z_2}\bT-\bS_i)^{-1}\\
&=&\underbrace{b_i(z_1)\sum\limits_{k<i}
\balpha_k^*\rE_i(\frac{1}{z_1}\bT-\bS_{ik})^{-1}\rE_i\bigg[(\frac{1}{z_2}\bT-\bS_i)^{-1}-
(\frac{1}{z_2}\bT-\bS_{ik})^{-1}\bigg](\frac{n-1}{n}b_i(z_1)\bI-\frac{1}{z_1}\bT)^{-1}\balpha_k}_{\bC_1}\\
&&\underbrace{-b_i(z_1)\rtr
\frac{1}{n}\sum\limits_{k<i}(\frac{n-1}{n}b_i(z_1)\bI-\frac{1}{z_1}\bT)^{-1}
\rE_i(\frac{1}{z_1}\bT-\bS_{ik})^{-1}\rE_i[(\frac{1}{z_2}\bT-\bS_i)^{-1}-(\frac{1}{z_2}\bT-\bS_{ik})^{-1}]}_{\bC_2}\\
&&+\underbrace{b_i(z_1)\rtr
\sum\limits_{k<i}(\frac{n-1}{n}b_i(z_1)\bI-\frac{1}{z_1}\bT)^{-1}
(\balpha_k\balpha_k^*-\frac{1}{n}\bI)rE_i(\frac{1}{z_1}\bT-\bS_{ik})^{-1}\rE_i(\frac{1}{z_2}\bT-\bS_{ik})^{-1}}_{\bC_3}\\
\end{array}
$$
where
$$
\begin{array}{lll}
\rE|\bC_2|&=&\rE\left|b_i(z_1)\rtr
\frac{1}{n}\sum\limits_{k<i}(\frac{n-1}{n}b_i(z_1)\bI-\frac{1}{z_1}\bT)^{-1}
\rE_i(\frac{1}{z_1}\bT-\bS_{ik})^{-1}rE_i[(\frac{1}{z_2}\bT-\bS_i)^{-1}-(\frac{1}{z_2}\bT-\bS_{ik})^{-1}]\right|\\
&=&\rE\bigg|b_i(z_1)\frac{1}{n}\sum\limits_{k<i}(\frac{n-1}{n}b_i(z_1)\bI-\frac{1}{z_1}\bT)^{-1}
\rE_i\beta_{k(i)}\balpha_k'(\frac{1}{z_1}\bT-\bS_{ik})^{-1}
(\frac{1}{z_1}\bT-\breve{\bS}_{ik})^{-1}(\frac{1}{z_1}\bT-\bS_{ik})^{-1}
\balpha_k\bigg|\\
&\leq& K\quad(\mbox{similar to the proof of (\ref{E1})})
\end{array}
$$
$$
\begin{array}{lll}
\rE|\bC_3|&=&\rE\left|b_i(z_1)\rtr
\sum\limits_{k<i}(\frac{n-1}{n}b_i(z_1)\bI-\frac{1}{z_1}\bT)^{-1}
(\balpha_k\balpha_k^*-\frac{1}{n}\bI)\rE_i(\frac{1}{z_1}\bT-\bS_{ik})^{-1}\rE_i(\frac{1}{z_2}\bT-\bS_{ik})^{-1}\right|\\
&\leq&\sum\limits_{k<i}K\rE\Bigg|
\balpha_k^*\rE_i(\frac{1}{z_1}\bT-\bS_{ik})^{-1}(\frac{1}{z_2}\bT-\breve{\bS}_{ik})^{-1}\balpha_k
-\frac{1}{n}\rtr\rE_i(\frac{1}{z_1}\bT-\bS_{ik})^{-1}(\frac{1}{z_2}\bT-\breve{\bS}_{ik})^{-1}\Bigg|\\
&\leq&\sum\limits_{k<i}KE^{\frac12}\Bigg|
\balpha_k^*\rE_i(\frac{1}{z_1}\bT-\bS_{ik})^{-1}(\frac{1}{z_2}\bT-\breve{\bS}_{ik})^{-1}\balpha_k
-\frac{1}{n}\rtr\rE_i(\frac{1}{z_1}\bT-\bS_{ik})^{-1}(\frac{1}{z_2}\bT-\breve{\bS}_{ik})^{-1}\Bigg|^2\\
&=&O(n^{\frac12})\quad\mbox{(by (1.15) of Bai and Silverstein
(2004))}
\end{array}
$$
$$
\begin{array}{lll}
\bC_1&=&b_i(z_1)\sum\limits_{k<i}
\balpha_k^*\rE_i(\frac{1}{z_1}\bT-\bS_{ik})^{-1}\rE_i[(\frac{1}{z_2}\bT-\bS_i)^{-1}-
(\frac{1}{z_2}\bT-\bS_{ik})^{-1}](\frac{n-1}{n}b_i(z_1)\bI-\frac{1}{z_1}\bT)^{-1}\balpha_k\\
 &=&b_i(z_1)\sum\limits_{k<i}
\rE_i\breve{\beta}_{k(i)}(z_2)\cdot\balpha_k^*(\frac{1}{z_1}\bT-\bS_{ik})^{-1}(\frac{1}{z_2}\bT-\breve{\bS}_{ik})^{-1}
\balpha_k\\
&&\cdot\balpha_k^*(\frac{1}{z_2}\bT-\breve{\bS}_{ik})^{-1}(\frac{n-1}{n}b_i(z_1)\bI-\frac{1}{z_1}\bT)^{-1}
\balpha_k\\
&=&b_i(z_1)b_i(z_2)\sum\limits_{k<i}
\rE_i\balpha_k^*(\frac{1}{z_1}\bT-\bS_{ik})^{-1}(\frac{1}{z_2}\bT-\breve{\bS}_{ik})^{-1}
\balpha_k\\
&&\cdot\balpha_k^*(\frac{1}{z_2}\bT-\breve{\bS}_{ik})^{-1}(\frac{n-1}{n}b_i(z_1)\bI-\frac{1}{z_1}\bT)^{-1}
\balpha_k+O_p(n^{1/2})\quad(\mbox{by (\ref{E2}) and (\ref{E3})})\\
&=&\frac{1}{n^2}b_i(z_1)b_i(z_2)\sum\limits_{k<i}
\rE_i\rtr(\frac{1}{z_1}\bT-\bS_{ik})^{-1}(\frac{1}{z_2}\bT-\breve{\bS}_{ik})^{-1}\\
 &&\cdot\rtr(\frac{1}{z_2}\bT-\breve{\bS}_{ik})^{-1}(\frac{n-1}{n}b_i(z_1)\bI-\frac{1}{z_1}\bT)^{-1}+O_p(n^{1/2})\\
&=&\frac{i-1}{n^2}b_i(z_1)b_i(z_2)
\rE_i\rtr(\frac{1}{z_1}\bT-\bS_{i})^{-1}(\frac{1}{z_2}\bT-\breve{\bS}_{i})^{-1}\\
 &&\cdot\rtr(\frac{1}{z_2}\bT-\breve{\bS}_{i})^{-1}(\frac{n-1}{n}b_i(z_1)\bI-\frac{1}{z_1}\bT)^{-1}+O_p(n^{1/2})\\
\end{array}
$$
That is,
$$
\begin{array}{lll}
&&\rtr[\rE_i(\frac{1}{z_1}\bT-\bS_{i})^{-1}](\frac{1}{z_2}\bT-\breve{\bS}_{i})^{-1}
\left[1+\frac{(i-1)}{n^2}b_i(z_1)b_i(z_2)
\cdot\rtr(\frac{1}{z_2}\bT-\breve{\bS}_{i})^{-1}(\frac{n-1}{n}b_i(z_1)\bI-\frac{1}{z_1}\bT)^{-1}\right]\\
&=&-\rtr\left[(\frac{n-1}{n}b_i(z_1)\bI-\frac{1}{z_1}\bT)^{-1}
\cdot\rE_i(\frac{1}{z_2}\bT-\bS_{i})^{-1}\right]+O_p(n^{1/2}).
\end{array}
$$
Then by (\ref{inverse}) we have
$$
\begin{array}{lll}
&&\rtr[\rE_i(\frac{1}{z_1}\bT-\bS_{i})^{-1}](\frac{1}{z_2}\bT-\breve{\bS}_{i})^{-1}\\
&&\left[1-\frac{(i-1)}{n^2}b_i(z_1)b_i(z_2)\rtr(\frac{n-1}{n}b_i(z_1)\bI-\frac{1}{z_1}\bT)^{-1}
\cdot(\frac{n-1}{n}b_i(z_2)\bI-\frac{1}{z_2}\bT)^{-1}\right]\\
&=&\rtr\left[(\frac{n-1}{n}b_i(z_1)\bI-\frac{1}{z_1}\bT)^{-1}
\cdot(\frac{n-1}{n}b_i(z_2)\bI-\frac{1}{z_2}\bT)^{-1}\right]+O_p(n^{1/2})
\end{array}
$$
because
$\rtr\bbA(z_1)(\frac{n-1}{n}b_i(z_1)\bI-\frac{1}{z_1}\bT)^{-1}=O_p(n^{1/2})$,
$\rtr\bbB(z_1)(\frac{n-1}{n}b_i(z_1)\bI-\frac{1}{z_1}\bT)^{-1}=O_p(n^{1/2})$
and
$\rtr\bbC(z_1)(\frac{n-1}{n}b_i(z_1)\bI-\frac{1}{z_1}\bT)^{-1}=O_p(n^{1/2})$.

By Lemma \ref{BaiLem91} and \ref{2}, we have
$$
|b_i(z)-b(z)|\leq Kn^{-1},\quad |b_i(z)-\rE\beta_i(z)|\leq
Kn^{-1/2},
$$
$$
\frac{1}{pz}\sum\limits_{i=1}^n\rE(-1+\beta_i(z))=\rE s_n(z),\quad
{\color{red}|\rE s_n(z)-s_n^0(z)|\leq Kn^{-1}}
$$
$$
E\beta_i(z)=y_nz\rE s_n(z)+1
$$
So we have
$$
\begin{array}{lll}
&&\rtr[\rE_i(\frac{1}{z_1}\bT-\bS_{i})^{-1}](\frac{1}{z_2}\bT-\breve{\bS}_{i})^{-1}\\
&&\left[1-\frac{(i-1)}{n^2}b(z_1)b(z_2)\rtr(b(z_1)\bI-\frac{1}{z_1}\bT)^{-1}
\cdot(b(z_2)\bI-\frac{1}{z_2}\bT)^{-1}\right]\\
&=&\rtr\left[b(z_1)\bI-\frac{1}{z_1}\bT)^{-1}
\cdot(b(z_2)\bI-\frac{1}{z_2}\bT)^{-1}\right]+O_p(n^{1/2})
\end{array}
$$
So we obtain
$$
\begin{array}{lll}
&&b(z_1)b(z_2)\rtr[\rE_i(\frac{1}{z_1}\bT-\bS_{i})^{-1}](\frac{1}{z_2}\bT-\breve{\bS}_{i})^{-1}
\left[1-\frac{(i-1)}{n}y_n\int\frac{b(z_1)b(z_2)}{(b(z_1)-\frac{1}{z_1}t)\cdot
(b(z_2)-\frac{1}{z_2}t)}dH_n(t)\right]\\
&=&p\int\frac{b(z_1)b(z_2)}{(b(z_1)-\frac{1}{z_1}t)\cdot
(b(z_2)-\frac{1}{z_2}t)}dH_n(t)+O_p(n^{1/2})
\end{array}
$$
That is,
$$
b(z_1)b(z_2)\rtr[\rE_i(\frac{1}{z_1}\bT-\bS_{i})^{-1}](\frac{1}{z_2}\bT-\breve{\bS}_{i})^{-1}=
\frac{p\int\frac{b(z_1)b(z_2)}{(b(z_1)-\frac{1}{z_1}t)\cdot
(b(z_2)-\frac{1}{z_2}t)}dH_n(t)}{1-\frac{(i-1)}{n}y_n\int\frac{b(z_1)b(z_2)}{(b(z_1)-\frac{1}{z_1}t)\cdot
(b(z_2)-\frac{1}{z_2}t)}dH_n(t)}
$$
Moreover, we have \bqn
&&\frac{b(z_1)b(z_2)}{n^2}\sum\limits_{i=1}^n\rtr[\rE_i(\frac{1}{z_1}\bT-\bS_{i})^{-1}](\frac{1}{z_2}\bT-\breve{\bS}_{i})^{-1}
\rightarrow\int\limits_0^1\frac{y\int\frac{b(z_1)b(z_2)}{(b(z_1)-\frac{1}{z_1}t)\cdot
(b(z_2)-\frac{1}{z_2}t)}dH(t)}{1-x\cdot
y\int\frac{b(z_1)b(z_2)}{(b(z_1)-\frac{1}{z_1}t)\cdot
(b(z_2)-\frac{1}{z_2}t)}dH(t)}dx\\
&=&a(z_1,z_2)\int\limits_0^1\frac{1}{1-x a(z_1,z_2)}dx=\int\limits_0^{a(z_1,z_2)}\frac{1}{1-z}dz
 \eqn
where
%$$
%a(z_1,z_2)=y\int\frac{b(z_1)b(z_2)}{(b(z_1)-\frac{1}{z_1}t)\cdot
%(b(z_2)-\frac{1}{z_2}t)}dH(t)
%$$
\bqn
a(z_1,z_2)&=&y\int\frac{b(z_1)b(z_2)}{(\frac{t}{z_1}-b(z_1))\cdot
(\frac{t}{z_2}-b(z_2))}dH(t)\\
&=&y\int\frac{(1+yz_1s(z_1))\cdot(1+yz_2s(z_2))}{(\frac{t}{z_1}-1-yz_1s(z_1))\cdot
(\frac{t}{z_2}-1-yz_2s(z_2))}dH(t)\\
&=&\frac{y}{z_1z_2}\int\frac{\underline{\tilde{m}}(\frac{1}{z_1})\cdot \underline{\tilde{m}}(\frac{1}{z_2})}{(\frac{t}{z_1}+\frac{\underline{\tilde{m}}(\frac{1}{z_1})}{z_1})\cdot
(\frac{t}{z_2}+\frac{\underline{\tilde{m}}(\frac{1}{z_2})}{z_2})}dH(t)\\
&=&y\int\frac{\underline{\tilde{m}}(\frac{1}{z_1})\cdot \underline{\tilde{m}}(\frac{1}{z_2})}{(t+\underline{\tilde{m}}(\frac{1}{z_1}))\cdot
(t+\underline{\tilde{m}}(\frac{1}{z_2}))}dH(t)\\
&=&\frac{\underline{\tilde{m}}(\frac{1}{z_1})\cdot \underline{\tilde{m}}(\frac{1}{z_2})}
{\underline{\tilde{m}}(\frac{1}{z_2})-\underline{\tilde{m}}(\frac{1}{z_1})}\cdot y\cdot\left(\int\frac{1}{t+\underline{\tilde{m}}(\frac{1}{z_1})}dH(t)-
\int\frac{1}{t+\underline{\tilde{m}}(\frac{1}{z_2})}dH(t)\right)\\
&=&\frac{\underline{\tilde{m}}(\frac{1}{z_1})\cdot \underline{\tilde{m}}(\frac{1}{z_2})}
{\underline{\tilde{m}}(\frac{1}{z_2})-\underline{\tilde{m}}(\frac{1}{z_1})}\left(
\frac{1}{z_1}-\frac{1}{z_2}+\frac{\underline{\tilde{m}}(\frac{1}{z_2})-\underline{\tilde{m}}(\frac{1}{z_1})}{\underline{\tilde{m}}(\frac{1}{z_1})\cdot\underline{\tilde{m}}(\frac{1}{z_2})}\right)\\
&=&1+\frac{\underline{\tilde{m}}(\frac{1}{z_1})\cdot \underline{\tilde{m}}(\frac{1}{z_2})}
{\underline{\tilde{m}}(\frac{1}{z_2})-\underline{\tilde{m}}(\frac{1}{z_1})}\left(
\frac{1}{z_1}-\frac{1}{z_2}\right)\\
\eqn
 with $\underline{\tilde{m}}(\frac{1}{z})\stackrel{\triangle}{=}-z(1+yzs(z))$ and
$$
\begin{array}{lll}
 \int\frac{1}{t+\underline{\tilde{m}}(\frac{1}{z})}dH(t)&=&
\frac{1}{z}\int\frac{1}{\frac{t}{z}+\frac1z\underline{\tilde{m}}(\frac{1}{z})}dH(t)\\
&=&\frac{1}{z}\int\frac{1}{\frac{t}{z}-(1+yzs(z))}dH(t)=\frac{\tilde{m}(z)}{z}\quad(\mbox{by}~(\ref{eq24})~\mbox{and}~(\ref{7}))\\
&=&\frac{s(z)}{1+yzs(z)}=
\frac1y\left(\frac1z+\frac{1}{\underline{\tilde{m}}(\frac{1}{z})}\right)\quad(\mbox{by}~(\ref{7})).
\end{array}
$$
That is, the covariance is
$$
\frac{\partial^2}{\partial\alpha_1\partial\alpha_2}\int\limits_0^{a(z_1,z_2)}\frac{1}{1-z}dz
=\frac{\partial}{\partial z_2}\left(\frac{\partial a(z_1,z_2)/\partial z_1}{1-a(z_1,z_2)}\right)
$$
\bqn
\partial a(z_1,z_2)/\partial z_1&=&\left(\frac{\underline{\tilde{m}}(\frac{1}{z_1})\cdot \underline{\tilde{m}}(\frac{1}{z_2})}
{\underline{\tilde{m}}(\frac{1}{z_2})-\underline{\tilde{m}}(\frac{1}{z_1})}\left(
\frac{1}{z_1}-\frac{1}{z_2}\right)\right)_{z_1}\\
&=&-\frac{(\underline{\tilde{m}}(\frac{1}{z_1}))'\underline{\tilde{m}}(\frac{1}{z_2})(\underline{\tilde{m}}(\frac{1}{z_2})-\underline{\tilde{m}}(\frac{1}{z_1}))+\underline{\tilde{m}}(\frac{1}{z_1}))'\underline{\tilde{m}}(\frac{1}{z_1}))\underline{\tilde{m}}(\frac{1}{z_2})}
{(\underline{\tilde{m}}(\frac{1}{z_2})-\underline{\tilde{m}}(\frac{1}{z_1}))^2}\left(
\frac{1}{z_1}-\frac{1}{z_2}\right)\\
&&+\frac{\underline{\tilde{m}}(\frac{1}{z_1})\cdot \underline{\tilde{m}}(\frac{1}{z_2})}
{\underline{\tilde{m}}(\frac{1}{z_2})-\underline{\tilde{m}}(\frac{1}{z_1})}\left(
\frac{-1}{z_1^2}\right)\\
&=&-\frac{(\underline{\tilde{m}}(\frac{1}{z_1}))'\underline{\tilde{m}}^2(\frac{1}{z_2})}
{(\underline{\tilde{m}}(\frac{1}{z_2})-\underline{\tilde{m}}(\frac{1}{z_1}))^2}\left(
\frac{1}{z_1}-\frac{1}{z_2}\right)+\frac{\underline{\tilde{m}}(\frac{1}{z_1})\cdot
\underline{\tilde{m}}(\frac{1}{z_2})}
{\underline{\tilde{m}}(\frac{1}{z_2})-\underline{\tilde{m}}(\frac{1}{z_1})}\left(
\frac{-1}{z_1^2}\right). \eqn
So we obtain \bqn &&\frac{\partial
a(z_1,z_2)/\partial
z_1}{1-a(z_1,z_2)}\\
&=&\left[-\frac{(\underline{\tilde{m}}(\frac{1}{z_1}))'\underline{\tilde{m}}^2(\frac{1}{z_2})}
{(\underline{\tilde{m}}(\frac{1}{z_2})-\underline{\tilde{m}}(\frac{1}{z_1}))^2}\left(
\frac{1}{z_1}-\frac{1}{z_2}\right)+\frac{\underline{\tilde{m}}(\frac{1}{z_1})\cdot
\underline{\tilde{m}}(\frac{1}{z_2})}
{\underline{\tilde{m}}(\frac{1}{z_2})-\underline{\tilde{m}}(\frac{1}{z_1})}\left(
\frac{-1}{z_1^2}\right)\right]\frac
{\underline{\tilde{m}}(\frac{1}{z_2})-\underline{\tilde{m}}(\frac{1}{z_1})}{\underline{\tilde{m}}(\frac{1}{z_1})\underline{\tilde{m}}(\frac{1}{z_2})\left(
\frac{1}{z_1}-\frac{1}{z_2}\right)}\\
&=&-\frac{(\underline{\tilde{m}}(\frac{1}{z_1}))'\underline{m}(\frac{1}{z_2})}{m(\frac{1}{z_1})(\underline{\tilde{m}}(\frac{1}{z_2})-\underline{\tilde{m}}(\frac{1}{z_1}))}
-\frac{1/z_1^2}{1/z_1-1/z_2} \eqn and \bqn
\frac{\partial}{\partial z_2}\left(\frac{\partial
a(z_1,z_2)/\partial z_1}{1-a(z_1,z_2)}\right)
&=&\frac{(\underline{\tilde{m}}(\frac{1}{z_1}))'(\underline{\tilde{m}}(\frac{1}{z_2}))'}{(\underline{\tilde{m}}(\frac{1}{z_2})-\underline{\tilde{m}}(\frac{1}{z_1}))^2}
-\frac{1}{z_1^2z_2^2}\frac{1}{(1/z_1-1/z_2)^2}\\
&=&\frac{(\underline{\tilde{m}}(\frac{1}{z_1}))'(\underline{\tilde{m}}(\frac{1}{z_2}))'}{(\underline{\tilde{m}}(\frac{1}{z_2})-\underline{\tilde{m}}(\frac{1}{z_1}))^2}
-\frac{1}{(z_1-z_2)^2}\\
&=&\frac{[z_1(1+yz_1s(z_1))]'[z_2(1+yz_2s(z_2))]'}{[z_1(1+yz_1s(z_1))-z_2(1+yz_2s(z_2))]^2}
-\frac{1}{(z_1-z_2)^2}. \eqn That is, by (\ref{var}) we have
\begin{eqnarray}
\mbox{Cov}(M(z_1),
M(z_2))=\kappa\left(\frac{[z_1(1+yz_1s(z_1))]'[z_2(1+yz_2s(z_2))]'}{[z_1(1+yz_1s(z_1))-z_2(1+yz_2s(z_2))]^2}
-\frac{1}{(z_1-z_2)^2}\right). \label{var_complement}
\end{eqnarray}

Then the proof of Theorem \ref{thm3} is completed. \eprf

\subsection{Tightness of $M_n^{1}(z)$}

\begin{thm}\label{thm4}
Under Assumptions 1$^*$, 2$^*$ and 3, the sequence of random
functions $M^{1}_n(z)$ is tight for
$z\in\mathcal{C}\cup\bar{\mathcal{C}}$.
\end{thm}
\proof We want to show that
%$$\sup\limits_{n;z_1,z_2\in\mathcal{C}_n}\frac{\rE|\hat{M}_n^1(z_1)-\hat{M}_n^1(z_2))|^2}{|z_1-z_2|^2}=$$
$$
\sup\limits_{n;z_1,z_2\in
\mathcal{C}_n}\frac{\rE|M_n^1(z_1)-M_n^1(z_2))|^2}{|z_1-z_2|^2}
$$
is finite. It is straightforward to verify that this will be true
if we can find a $K>0$ for which
$$
\sup\limits_{n;z_1,z_2\in
\mathcal{C}_n}\frac{\rE|M_n^1(z_1)-M_n^1(z_2))|^2}{|z_1-z_2|^2}\leq
K.
$$
$$
M_n^1(z_1)-M_n^1(z_2)=p(s_n(z_1)-s_n(z_2))-p\cdot
\rE(s_n(z_1)-s_n(z_2))
$$
$$
s_n(z_1)-s_n(z_2)=\frac1p\rtr\left[(\bS^{-1}\bT-z_1\bI)^{-1}-(\bS^{-1}\bT-z_2\bI)^{-1}\right]=\frac1p(z_1-z_2)\rtr\bD^{-1}(z_1)\bD^{-1}(z_2)
$$
where $\bD(z)=(\bS^{-1}\bT-z\bI)^{-1}$. We have
$$
\begin{array}{lll}
\displaystyle{\frac{M_n^1(z_1)-M_n^1(z_2)}{z_1-z_2}}&=&\displaystyle{p\cdot\frac{s_n(z_1)-s_n(z_2)-\rE(s_n(z_1)-s_n(z_2))}{z_1-z_2}}\\
&=&\sum\limits_{i=1}^n(\rE_i-\rE_{i-1})\rtr\bD^{-1}(z_1)\bD^{-1}(z_2)\\
&=&\sum\limits_{i=1}^n(\rE_i-E_{i-1})\rtr(\bS^{-1}\bT-z_1\bI)^{-1}(\bS^{-1}\bT-z_2\bI)^{-1}\\
\end{array}
$$
and
$$
\begin{array}{lll}
(\bS^{-1}\bT-z\bI)^{-1}&=&\bS(\bT-z\bS)^{-1}\\
&=&\bS(\bT-z\bS_i-z\bm{\alpha}_i\bm{\alpha}_i')^{-1}\\
&=&(\bS_i+\bm{\alpha}_i\bm{\alpha}_i')\left(\bD_i^{-1}+\frac{z\cdot\bD_i^{-1}\bm{\alpha}_i\bm{\alpha}_i'\bD_i^{-1}}{1-z\cdot\bm{\alpha}_i'\bD_i^{-1}
\bm{\alpha}_i}\right)\\
&=&\bS_i\bD_i^{-1}+\frac{z\cdot\bS_i\bD_i^{-1}\bm{\alpha}_i\bm{\alpha}_i'\bD_i^{-1}}{1-z\cdot\bm{\alpha}_i'\bD_i^{-1}\bm{\alpha}_i}+
\bm{\alpha}_i\bm{\alpha}_i'\bD_i^{-1}+\frac{z\cdot\bm{\alpha}_i\bm{\alpha}_i'\bD_i^{-1}\cdot\bm{\alpha}_i'\bD_i^{-1}\bm{\alpha}_i}
{1-z\cdot\bm{\alpha}_i'\bD_i^{-1}\bm{\alpha}_i}\\
&=&\bS_i\bD_i^{-1}+\frac{z\cdot\bS_i\bD_i^{-1}\bm{\alpha}_i\bm{\alpha}_i'\bD_i^{-1}}{1-z\cdot\bm{\alpha}_i'\bD_i^{-1}\bm{\alpha}_i}+
\frac{\bm{\alpha}_i\bm{\alpha}_i'\bD_i^{-1}}
{1-z\cdot\bm{\alpha}_i'\bD_i^{-1}\bm{\alpha}_i}\\
&=&(\bS_i^{-1}\bT-z\bI)^{-1}+
\frac{z\cdot\bS_i\bD_i^{-1}\bm{\alpha}_i\bm{\alpha}_i'\bD_i^{-1}}{1-z\cdot\bm{\alpha}_i'\bD_i^{-1}\bm{\alpha}_i}+
\frac{\bm{\alpha}_i\bm{\alpha}_i'\bD_i^{-1}}
{1-z\cdot\bm{\alpha}_i'\bD_i^{-1}\bm{\alpha}_i}\\
&=&(\bS_i^{-1}\bT-z\bI)^{-1}+
\frac{(z\cdot\bS_i\bD_i^{-1}+\bI)\bm{\alpha}_i\bm{\alpha}_i'\bD_i^{-1}}{1-z\cdot\bm{\alpha}_i'\bD_i^{-1}\bm{\alpha}_i}\\
\end{array}
$$
where $\bS_i\bD_i^{-1}=(\bS_i^{-1}\bT-z\bI)^{-1}$. That is,
$$
(\bS^{-1}\bT-z\bI)^{-1}-(\bS_i^{-1}\bT-z\bI)^{-1}=
\frac{(z\cdot\bS_i\bD_i^{-1}+\bI)\bm{\alpha}_i\bm{\alpha}_i'\bD_i^{-1}}{1-z\cdot\bm{\alpha}_i'\bD_i^{-1}\bm{\alpha}_i}=
(z\bS_i\bD_i^{-1}+\bI)\bm{\alpha}_i\bm{\alpha}_i'\bD_i^{-1}\beta_i(z).
$$
$$
\begin{array}{lll}
&&\rtr(\bS^{-1}\bT-z_1\bI)^{-1}(\bS^{-1}\bT-z_2\bI)^{-1}-\rtr(\bS_i^{-1}\bT-z_1\bI)^{-1}(\bS_i^{-1}\bT-z_2\bI)^{-1}\\
&=&\rtr[(\bS^{-1}\bT-z_1\bI)^{-1}-(\bS_i^{-1}\bT-z_1\bI)^{-1}][(\bS^{-1}\bT-z_2\bI)^{-1}-(\bS_i^{-1}\bT-z_2\bI)^{-1}]\\
&&+\rtr[(\bS^{-1}\bT-z_1\bI)^{-1}-(\bS_i^{-1}\bT-z_1\bI)^{-1}](\bS_i^{-1}\bT-z_2\bI)^{-1}\\
&&+\rtr(\bS_i^{-1}\bT-z_1\bI)^{-1}[(\bS^{-1}\bT-z_2\bI)^{-1}-(\bS_i^{-1}\bT-z_2\bI)^{-1}]\\
&=&\bm{\alpha}_i'\bD_i^{-1}
(z_2\bS_i\bD_i^{-1}+\bI)\bm{\alpha}_i\cdot\bm{\alpha}_i'\bD_i^{-1}(z_1\bS_i\bD_i^{-1}+\bI)\bm{\alpha}_i\cdot\beta_i(z_1)\beta_i(z_2)\\
&&+\beta_i(z_1)\cdot\bm{\alpha}_i'\bD_i^{-1}(\bS_i^{-1}\bT-z_2\bI)^{-1}(z_1\bS_i\bD_i^{-1}+\bI)\bm{\alpha}_i\\
&&+\beta_i(z_2)\cdot\bm{\alpha}_i'\bD_i^{-1}(\bS_i^{-1}\bT-z_1\bI)^{-1}(z_2\bS_i\bD_i^{-1}+\bI)\bm{\alpha}_i\\
&=&(\bm{\alpha}_i'\bD_i^{-1}
\bS_i^{-1}\bT{\bf F}_i^{-1}(z_2)\bm{\alpha}_i)^2\cdot\beta_i(z_1)\beta_i(z_2)+\beta_i(z_1)\cdot\bm{\alpha}_i'\bD_i^{-1}{\bf F}_i^{-1}(z_2)\bS_i^{-1}\bT{\bf F}_i^{-1}(z_1)\bm{\alpha}_i\\
&&+\beta_i(z_2)\cdot\bm{\alpha}_i'\bD_i^{-1}{\bf F}_i^{-1}(z_1)\bS_i^{-1}\bT{\bf F}_i^{-1}(z_2)\bm{\alpha}_i\\
\end{array}
$$
where ${\bf F}^{-1}_i(z)=(\bS_i^{-1}\bT-z\bI)^{-1}$ and
$(\bS_i^{-1}\bT-z\bI)^{-1}=-\frac1z\bI+\frac1z\bS_i^{-1}\bT(\bS_i^{-1}\bT-z\bI)^{-1}$.

\bqa
\frac{M_n^1(z_1)-M_n^1(z_2)}{z_1-z_2}
&=&\sum\limits_{i=1}^n(\rE_i-\rE_{i-1})\rtr(\bS^{-1}\bT-z_1\bI)^{-1}(\bS^{-1}\bT-z_2\bI)^{-1}\label{one}\\
&=&\sum\limits_{i=1}^n(\rE_i-\rE_{i-1})(\bm{\alpha}_i'\bD_i^{-1}
\bS_i^{-1}\bT{\bf
F}_i^{-1}(z_2)\bm{\alpha}_i)^2\cdot\beta_i(z_1)\beta_i(z_2)\label{one1}\\
&&+\sum\limits_{i=1}^n(\rE_i-\rE_{i-1})\beta_i(z_1)\cdot\bm{\alpha}_i'\bD_i^{-1}{\bf
F}_i^{-1}(z_2)\bS_i^{-1}\bT{\bf F}_i^{-1}(z_1)\bm{\alpha}_i\label{one2}\\
&&+\sum\limits_{i=1}^n(\rE_i-\rE_{i-1})\beta_i(z_2)\cdot\bm{\alpha}_i'\bD_i^{-1}{\bf
F}_i^{-1}(z_1)\bS_i^{-1}\bT{\bf F}_i^{-1}(z_2)\bm{\alpha}_i\label{one3}
\eqa
Our goal is to show that the absolute second moment of (\ref{one})
is bounded. We begin with (\ref{one2}). We have
$$
\begin{array}{lll}
&&\sum\limits_{i=1}^n(\rE_i-\rE_{i-1})\beta_i(z_1)\cdot\bm{\alpha}_i'\bD_i^{-1}{\bf
F}_i^{-1}(z_2)\bS_i^{-1}\bT{\bf F}_i^{-1}(z_1)\bm{\alpha}_i\\
&=&\sum\limits_{i=1}^n(\rE_i-\rE_{i-1})b_i(z_1)\cdot\bm{\alpha}_i'\bD_i^{-1}{\bf
F}_i^{-1}(z_2)\bS_i^{-1}\bT{\bf F}_i^{-1}(z_1)\bm{\alpha}_i\\
&&-\sum\limits_{i=1}^n(\rE_i-\rE_{i-1})\beta_i(z_1)b_i(z_1)\varepsilon_i(z_1)\cdot\bm{\alpha}_i'\bD_i^{-1}{\bf
F}_i^{-1}(z_2)\bS_i^{-1}\bT{\bf F}_i^{-1}(z_1)\bm{\alpha}_i\\
&=&W_1-W_2
\end{array}
$$
where
$\varepsilon_i(z)=\bm{\alpha}_i'(\frac1z\bT-\bS_i)^{-1}\bm{\alpha}_i-
\frac1n\rE\rtr(\frac1z\bT-\bS_i)^{-1}$ and
$$
\begin{array}{lll}
\rE|W_1|^2&=&\rE|\sum\limits_{i=1}^n(\rE_i-E_{i-1})b_i(z_1)\cdot\bm{\alpha}_i'\bD_i^{-1}{\bf
F}_i^{-1}(z_2)\bS_i^{-1}\bT{\bf F}_i^{-1}(z_1)\bm{\alpha}_i|^2\\
        &=&\sum\limits_{i=1}^n\rE|(\rE_i-E_{i-1})b_i(z_1)\cdot\bm{\alpha}_i'\bD_i^{-1}{\bf
F}_i^{-1}(z_2)\bS_i^{-1}\bT{\bf F}_i^{-1}(z_1)\bm{\alpha}_i|^2\\
        &=&\sum\limits_{i=1}^nb_i^2(z_1)\rE|\rE_i\bm{\alpha}_i'\bD_i^{-1}{\bf
F}_i^{-1}(z_2)\bS_i^{-1}\bT{\bf
F}_i^{-1}(z_1)\bm{\alpha}_i-\frac1n\rtr\bD_i^{-1}{\bf
F}_i^{-1}(z_2)\bS_i^{-1}\bT{\bf F}_i^{-1}(z_1)|^2\\
       &\leq& K
\end{array}
$$
by Lemma \ref{2}. Moreover, we have
$$
\begin{array}{lll}
\rE|W_2|^2&=&\sum\limits_{i=1}^n\rE|(\rE_i-\rE_{i-1})\beta_i(z_1)b_i(z_1)\varepsilon_i(z_1)\cdot\bm{\alpha}_i'\bD_i^{-1}{\bf
F}_i^{-1}(z_2)\bS_i^{-1}\bT{\bf F}_i^{-1}(z_1)\bm{\alpha}_i|^2\\
&=&n\rE|(\rE_i-E_{i-1})\beta_i(z_1)b_i(z_1)\varepsilon_i(z_1)\cdot\bm{\alpha}_i'\bD_i^{-1}{\bf
F}_i^{-1}(z_2)\bS_i^{-1}\bT{\bf F}_i^{-1}(z_1)\bm{\alpha}_i|^2\\
&\leq&2nb^2_i(z_1)\cdot\rE\beta_i^2(z_1)\varepsilon^2_i(z_1)\cdot\left(\bm{\alpha}_i'\bD_i^{-1}{\bf
F}_i^{-1}(z_2)\bS_i^{-1}\bT{\bf F}_i^{-1}(z_1)\bm{\alpha}_i\right)^2\leq K\\
\end{array}
$$
where
$$
\begin{array}{lll}
|\beta_i(z)|&=&\left|1-\bm{\alpha}_i'(\frac1z\bT-\bS_i)^{-1}\bm{\alpha}_i\right|
\leq 1+|z\bm{\alpha}_i\bS_i^{-1}(\bT\bS_i^{-1}-z\bI)^{-1}\bm{\alpha}_i|\\
            &\leq&1+|z|\cdot|\bm{\alpha}_i|^2\|\bS_i^{-1}(\bT\bS_i^{-1}-z\bI)^{-1}\|\\
            &\leq&1+K|\bm{\alpha}_i|^2+n^5I(\|\bS_i\|\geq \eta_r\mbox{or}\lambda_{\min}^{\bS_i}\leq\eta_l)
\end{array}
$$
and
$$
\bm{\alpha}_i'\bD_i^{-1}{\bf
F}_i^{-1}(z_2)\bS_i^{-1}\bT{\bf F}_i^{-1}(z_1)\bm{\alpha}_i\leq K|\bm{\alpha}_i|^2+n^tI(\|\bS_i\|\geq \eta_r\mbox{or}\lambda_{\min}^{\bS_i}\leq\eta_l).
$$
Similarly, it can be obtained that the second moment of (\ref{one3}) is uniformly finite.
Now we begin (\ref{one1}). We have
\bqn
&&\sum\limits_{i=1}^n(\rE_i-E_{i-1})(\bm{\alpha}_i'\bD_i^{-1}
\bS_i^{-1}\bT{\bf F}_i^{-1}(z_2)\bm{\alpha}_i)^2\cdot\beta_i(z_1)\beta_i(z_2)\\
&=&\sum\limits_{i=1}^n(\rE_i-E_{i-1})\bigg(\bigg[\bm{\alpha}_i'\bD_i^{-1}
\bS_i^{-1}\bT{\bf F}_i^{-1}(z_2)\bm{\alpha}_i\bigg]^2-\bigg[\frac1n\rtr(\bD_i^{-1}
\bS_i^{-1}\bT{\bf F}_i^{-1}(z_2)\bigg]^2\bigg)\cdot b_i(z_1)b_i(z_2)\\
&&-\sum\limits_{i=1}^n(\rE_i-E_{i-1})(\bm{\alpha}_i'\bD_i^{-1}
\bS_i^{-1}\bT{\bf F}_i^{-1}(z_2)\bm{\alpha}_i)^2\cdot \beta_i(z_1)\beta_i(z_2)b_i(z_2)\varepsilon_i(z_2)\\
&&-\sum\limits_{i=1}^n(\rE_i-\rE_{i-1})(\bm{\alpha}_i'\bD_i^{-1}
\bS_i^{-1}\bT{\bf F}_i^{-1}(z_2)\bm{\alpha}_i)^2\cdot b_i(z_1)b_i(z_2)\beta_i(z_1)\varepsilon_i(z_1)\\
&=&Z_1-Z_2-Z_3. \eqn By similar methods to $W_2$, we obtain that
the second moments of $Z_2$ and $Z_3$ are uniformly finite. Now we
begin $Z_1$. We have \bqn
&&\rE|Z_1|^2\\
&\leq& K\sum\limits_{i=1}^n\rE\left|\bigg[\bm{\alpha}_i'\bD_i^{-1}
\bS_i^{-1}\bT{\bf F}_i^{-1}(z_2)\bm{\alpha}_i\bigg]^2-\bigg[\frac1n\rtr(\bD_i^{-1}
\bS_i^{-1}\bT{\bf F}_i^{-1}(z_2))\bigg]^2\right|^2\\
&\leq&2K\sum\limits_{i=1}^n\rE\left|\bm{\alpha}_i'\bD_i^{-1}
\bS_i^{-1}\bT{\bf F}_i^{-1}(z_2)\bm{\alpha}_i-\frac1n\rtr(\bD_i^{-1}
\bS_i^{-1}\bT{\bf F}_i^{-1}(z_2))\right|^4\\
&&+\frac{K}{n}\sum\limits_{i=1}^n\rE\left|\bigg(\bm{\alpha}_i'\bD_i^{-1}
\bS_i^{-1}\bT{\bf
F}_i^{-1}(z_2)\bm{\alpha}_i-\frac1n\rtr(\bD_i^{-1}
\bS_i^{-1}\bT{\bf F}_i^{-1}(z_2))\bigg)\cdot|\bD_i^{-1}
\bS_i^{-1}\bT{\bf F}_i^{-1}(z_2)|\right|^2\\
&\leq& K.
\eqn
 Then $\sup\limits_{n;z_1,z_2\in
\mathcal{C}^{+}}\frac{\rE|M_n^1(z_1)-M_n^1(z_2))|^2}{|z_1-z_2|^2}\leq K$.

So the proof of Theorem \ref{thm4} is completed. \eprf

\subsection{Uniform convergence of $M_n^2(z)=p(\rE s_n(z)-s_n^0(z))$
for $z\in{\cal{C}}$}

\begin{thm}\label{thm5}
We have
$$
\sup\limits_{z\in\mathcal{C}_n}\left|M_n^2(z)-\frac{1}{z^2}\cdot \frac{y\int
\frac{t(1+yzs(z))^3dH(t)}{(t/z-1-yzs(z))^3}}{\left(
1-y\int\frac{(1+yzs(z))^2 dH(t)}{(t/z-1-yzs(z))^2}\right)^2}\right|\rightarrow0\quad\mbox{as}~n\rightarrow+\infty
$$
where $M_n^2(z)$ converges uniformly to the limit
$$\frac{1}{z^2}\cdot \frac{y\int
\frac{t(1+yzs(z))^3dH(t)}{(t/z-1-yzs(z))^3}}{\left(
1-y\int\frac{(1+yzs(z))^2
dH(t)}{(t/z-1-yzs(z))^2}\right)^2}=\frac12\frac{d\log\left(1-y\int\frac{(1+yzs(z))^2
dH(t)}{(t/z-1-yzs(z))^2}\right)}{dz}.$$
\end{thm}

\proof
We have
$$
\frac1p\rtr({\bf S}^{-1}{\bf T}-z{\bf I})^{-1}=-\frac1z-\frac{1}{z^2}\cdot\frac1p\rtr\left({\bf S}-\frac{{\bf T}}{z}\right)^{-1}{\bf T}
$$
and
$$
s_n(z)=\frac{1}{pz}\sum\limits_{i=1}^n\frac{\balpha_i'(\frac{1}{z}\bT-\bS_i)^{-1}\balpha_i}
         {1-\balpha_i'(\frac{1}{z}\bT-\bS_i)^{-1}\balpha_i}
=-\frac{1}{yz}+\frac{1}{yz}\cdot\frac1n\sum\limits_{i=1}^n\frac{1}
         {1-\balpha_i'(\frac{1}{z}\bT-\bS_i)^{-1}\balpha_i}.
$$
By Lemma \ref{lem4}, we have
$s(z)=-\frac1z-\frac{1}{z^2}\cdot\tilde{s}(z)$
where
$$
\underline{\tilde{s}}(z)=\frac{-z}{1-y\tilde{m}(z)}=\frac{-z}{1-y\int\frac{1}{\frac{t}{z}-\frac{1}{1-y\tilde{m}(z)}}dH(t)}
=\frac{1}{-\frac1z+y\int\frac{1}{t+\underline{\tilde{s}}(z)}dH(t)},
$$
$\tilde{m}(z)$ is the limit of
$\frac1p\rtr\left(\frac1z\bbT-\bbS\right)^{-1}$, $\tilde{s}(z)$ is
the limit of $\frac1p\rtr\left({\bf S}-\frac{{\bf
T}}{z}\right)^{-1}{\bf T}$ and
$\underline{\tilde{s}}(z)=-(1-y)z+y\tilde{s}(z).$ We have
$$
\rE\underline{\tilde{s}}_n(z)=\frac{1}{-\frac1z+y\int\frac{1}{t+\rE\underline{\tilde{s}}_n(z)}dH_n(t)-R_n}
$$
where
$$
R_n=-\frac1z+y\int\frac{1}{t+\rE\underline{\tilde{s}}_n(z)}dH_n(t)-\frac{1}{\rE\underline{\tilde{s}}_n(z)}.
$$
Moreover, let
$$
\underline{\tilde{s}}_n^0(z)=\frac{1}{-\frac1z+y\int\frac{1}{t+\underline{\tilde{s}}_n^0(z)}dH_n(t)}
$$
and
$$
s_n^0(z)=-\frac1z-\frac{1}{z^2}\cdot\tilde{s}_n^0(z), \quad
\underline{\tilde{s}}_n^0(z)=-(1-y)z+y\tilde{s}_n^0(z)
$$
where $s_n^0(z)$ is the Stieltjes transform of $F^{y_n,H_n}$. So
we obtain
$$
\rE\underline{\tilde{s}}_n(z)-\underline{\tilde{s}}_n^0(z)=
\frac{(\rE\underline{\tilde{s}}_n(z)-\underline{\tilde{s}}_n^0(z))y\int\frac{1}
{(t+\rE\underline{\tilde{s}}_n(z))(t+\underline{\tilde{s}}_n^0(z))}dH_n(t)}
{\left(-\frac1z+y\int\frac{1}{t+\rE\underline{\tilde{s}}_n(z)}dH_n(t)-R_n\right)
\left(-\frac1z+y\int\frac{1}{t+\underline{\tilde{s}}_n^0(z)}dH_n(t)\right)}
+\rE\underline{\tilde{s}}_n(z)\underline{\tilde{s}}_n^0(z)R_n
$$
That is,
$$
n(\rE\underline{\tilde{s}}_n(z)-\underline{\tilde{s}}_n^0(z))=
\frac{\rE\underline{\tilde{s}}_n(z)\underline{\tilde{s}}_n^0(z)}
{1-\frac{y\int\frac{1}
{(t+\rE\underline{\tilde{s}}_n(z))(t+\underline{\tilde{s}}_n^0(z))}dH_n(t)}
{\left(-\frac1z+y\int\frac{1}{t+\rE\underline{\tilde{s}}_n(z)}dH_n(t)-R_n\right)
\left(-\frac1z+y\int\frac{1}{t+\underline{\tilde{s}}_n^0(z)}dH_n(t)\right)}}\cdot nR_n
$$
$$
R_n=-\frac{1}{\rE\underline{\tilde{s}}_n(z)}\cdot\left(\frac{y}{z}\rE\tilde{s}_n(z)+y\int\frac{tdH_n(t)}
{t+\rE\underline{\tilde{s}}_n(z)}\right)=
-\frac{1}{\rE\underline{\tilde{s}}_n(z)}\cdot\frac{y}{z}\left(\rE\tilde{s}_n(z)+\int\frac{tdH_n(t)}
{t/z+\rE\underline{\tilde{s}}_n(z)/z}\right)
$$

$$
\bigg(\frac1z\bT-\bS\bigg)=\bigg(\frac1z\bT+\frac{\rE\underline{\tilde{s}}_n(z)}{z}\bI\bigg)
-\frac{\rE\underline{\tilde{s}}_n(z)}{z}{\bf I}-
\sum\limits_{i=1}^n\bm{\alpha}_i\bm{\alpha}_i'
$$
So we obtain
\bqn
&&\bigg(\frac1z\bT-\bS\bigg)^{-1}{\bf T}\\
&=&\bigg(\frac1z\bT-\rE\beta_i(z)\bI\bigg)^{-1}{\bf T}
+\bigg(\frac1z\bT-\rE\beta_i(z) \bI\bigg)^{-1}\bigg(\sum\limits_{i=1}^n\bm{\alpha}_i\bm{\alpha}_i'-E\beta_i(z){\bf T}\bigg)\bigg(\frac1z\bT-\bS\bigg)^{-1}{\bf T}\\
&=&\bigg(\frac1z\bT-\rE\beta_i(z)\bI\bigg)^{-1}{\bf T}
-\rE\beta_i(z)\cdot\bigg(\frac1z\bT-\rE\beta_i(z) \bI\bigg)^{-1}\bigg(\frac1z\bT-\bS\bigg)^{-1}{\bf T}\\
&&+\sum\limits_{i=1}^n\bigg(\frac1z\bT-\rE\beta_i(z) \bI\bigg)^{-1}\bm{\alpha}_i\bm{\alpha}_i'\bigg(\frac1z\bT-\bS_i\bigg)^{-1}{\bf T}\beta_i(z)\\
\eqn
where $\beta_i(z)=\frac{1}{1-\bm{\alpha}_i'(\frac{\bT}{z}-\bS_i)^{-1}\bm{\alpha}_i}$.
Taking expected values and trace on both sides and dividing by $p$, we get
\bqn
&&\frac1p\rE\rtr\bigg(\frac1z\bT-\bS\bigg)^{-1}{\bf T}\\
&=&\frac1p\rtr\bigg(\frac1z\bT-\rE\beta_1(z)\cdot
\bI\bigg)^{-1}{\bf T}-E\beta_1(z)\cdot\frac1p\rE\rtr\bigg(\frac1z\bT-E\beta_1(z)\cdot \bI\bigg)^{-1}\bigg(\frac1z\bT-\bS\bigg)^{-1}{\bf T}\\
&&+\frac1p\rE\rtr\sum\limits_{i=1}^n\bigg(\frac1z\bT-\rE\beta_1(z)\cdot \bI\bigg)^{-1}\bm{\alpha}_i\bm{\alpha}_i'\bigg(\frac1z\bT-\bS_i\bigg)^{-1}{\bf T}\beta_i(z)\\
&=&\frac1p\rtr\bigg(\frac1z\bT-\rE\beta_1(z)\cdot
\bI\bigg)^{-1}{\bf T}-\rE\beta_1(z)\cdot\frac1p\rE\rtr\bigg(\frac1z\bT-\rE\beta_1(z)\cdot \bI\bigg)^{-1}\bigg(\frac1z\bT-\bS\bigg)^{-1}{\bf T}\\
&&+\frac{1}{y_n}{\color{red}\rE\beta_1(z)}\bm{\alpha}_1'\bigg(\frac1z\bT-\bS_1\bigg)^{-1}{\bf T}\bigg(\frac1z\bT-\rE\beta_1(z) \bI\bigg)^{-1}\bm{\alpha}_1\\
&=&\frac1p\rtr\bigg(\frac1z\bT-\rE\beta_1(z)\cdot
\bI\bigg)^{-1}{\bf T}
+\frac{1}{y_n}\rE\beta_1(z)\bigg[\bm{\alpha}_1'\bigg(\frac1z\bT-\bS_1\bigg)^{-1}{\bf T}\bigg(\frac1z\bT-\rE\beta_1(z) \bI\bigg)^{-1}\bm{\alpha}_1\\
&&-\frac1n\rE\rtr\bigg(\frac1z\bT-\bS_1\bigg)^{-1}{\bf T}\bigg(\frac1z\bT-\rE\beta_1(z)\bI\bigg)^{-1}\bigg]\\
&&+\Bigg[-\frac{1}{y_n}\rE\beta_1(z)\frac1n\rE\rtr\bigg(\frac1z\bT-\rE\beta_1(z)\cdot \bI\bigg)^{-1}\bigg(\frac1z\bT-\bS\bigg)^{-1}{\bf T}\\
&&+\frac{1}{y_n}\rE\beta_1(z)\frac1n\rE\rtr\bigg(\frac1z\bT-\rE\beta_1(z)\cdot \bI\bigg)^{-1}\bigg(\frac1z\bT-\bS_1\bigg)^{-1}{\bf T}\Bigg]+o(\frac{1}{n})\\
&=&\frac1p\rtr\bigg(\frac1z\bT-\rE\beta_1(z)\cdot
\bI\bigg)^{-1}{\bf T}
+\frac{1}{y_n}\rE\beta_1(z)\bigg[\bm{\alpha}_1'\bigg(\frac1z\bT-\bS_1\bigg)^{-1}{\bf T}\bigg(\frac1z\bT-\rE\beta_1(z) \bI\bigg)^{-1}\bm{\alpha}_1\\
&&-\frac1n\rE\rtr\bigg(\frac1z\bT-\bS_1\bigg)^{-1}{\bf T}\bigg(\frac1z\bT-\rE\beta_1(z)\bI\bigg)^{-1}\bigg]\\
&&-\frac{(\rE\beta_1(z))^2}{y}\frac{1}{n^2}\rE\rtr\bigg(\frac1z\bT-\bS_1\bigg)^{-1}
\bigg(\frac1z\bT-\rE\beta_1(z)\bI\bigg)^{-1}
\bigg(\frac1z\bT-\bS_1\bigg)^{-1}{\bf T}+o(\frac{1}{n})\\
\eqn

So we obtain \bqa
\frac{-z\rE\underline{\tilde{s}}_n(z)}{y}\cdot R_n&=&\int\frac{tdH_n(t)}{\frac{t}{z}-\rE\beta_1(z)}-\rE\tilde{s}_n(z) \nonumber\\
   &=&\frac{1}{p}\rtr\bigg(\frac1z\bT-\rE\beta_1(z)\cdot \bI\bigg)^{-1}{\bf T}-
     \frac1p\rE\rtr\bigg(\frac1z\bT-\bS\bigg)^{-1}{\bf T} \nonumber\\
   &=&-\frac{1}{y_n}\rE\beta_1(z)\bigg[\bm{\alpha}_1'\bigg(\frac1z\bT-\bS_1\bigg)^{-1}{\bf T}
   \bigg(\frac1z\bT-\rE\beta_1(z) \bI\bigg)^{-1}\bm{\alpha}_1\nonumber\\
&&-\frac1n\rE\rtr\bigg(\frac1z\bT-\bS_1\bigg)^{-1}{\bf T}\bigg(\frac1z\bT-\rE\beta_1(z)\bI\bigg)^{-1}\bigg]\nonumber\\
&&+\frac{(\rE\beta_1(z))^2}{y}\frac{1}{n^2}\rE\rtr\bigg(\frac1z\bT-\bS_1\bigg)^{-2}{\bf
T}\bigg(\frac1z\bT-\rE\beta_1(z)\bI\bigg)^{-1}
+o(\frac{1}{n})\nonumber\\
&=&-\frac{1}{y_n}\rE\bar{\beta}^2_1(z)\bigg[\bm{\alpha}_1'\bigg(\frac1z\bT-\bS_1\bigg)^{-1}{\bf
T}
\bigg(\frac1z\bT-\rE\beta_1(z) \bI\bigg)^{-1}\bm{\alpha}_1\nonumber\\
&&-\frac1n\rE\rtr\bigg(\frac1z\bT-\bS_1\bigg)^{-1}{\bf
T}\bigg(\frac1z\bT-\rE\beta_1(z)\bI\bigg)^{-1}\bigg]
\varepsilon_j\nonumber\\
&&+\frac{(\rE\beta_1(z))^2}{y}\frac{1}{n^2}\rE\rtr\bigg(\frac1z\bT-\bS_1\bigg)^{-2}{\bf
T}\bigg(\frac1z\bT-\rE\beta_1(z)\bI\bigg)^{-1}
+o(\frac{1}{n})\nonumber\\
&=&-\frac{(\rE\beta_1(z))^2}{y}\cdot\rE\bigg[\bm{\alpha}_1'\bigg(\frac1z\bT-\bS_1\bigg)^{-1}{\bf
T}
\bigg(\frac1z\bT-\rE\beta_1(z) \bI\bigg)^{-1}\bm{\alpha}_1\nonumber\\
&&-\frac1n\rE\rtr\bigg(\frac1z\bT-\bS_1\bigg)^{-1}{\bf
T}\bigg(\frac1z\bT-\rE\beta_1(z)\bI\bigg)^{-1}\bigg]
\varepsilon_j\nonumber\\
&&+\frac{(\rE\beta_1(z))^2}{y}\frac{1}{n^2}\rE\rtr\bigg(\frac1z\bT-\bS_1\bigg)^{-2}{\bf
T}\bigg(\frac1z\bT-\rE\beta_1(z)\bI\bigg)^{-1}
+o(\frac{1}{n})\nonumber\\
&=&-\frac{(\rE\beta_1(z))^2}{y}\cdot\rE\bigg[\bm{\alpha}_1'\bigg(\frac1z\bT-\bS_1\bigg)^{-1}{\bf
T}
\bigg(\frac1z\bT-\rE\beta_1(z) \bI\bigg)^{-1}\bm{\alpha}_1\nonumber\\
&&-\frac1n\rE\rtr\bigg(\frac1z\bT-\bS_1\bigg)^{-1}{\bf
T}\bigg(\frac1z\bT-\rE\beta_1(z)\bI\bigg)^{-1}\bigg]
\varepsilon_j\nonumber\\
&&+\frac{(\rE\beta_1(z))^2}{y\cdot
n^2}\rE\rtr\bigg(\frac1z\bT-\bS_1\bigg)^{-2}{\bf T}
\bigg(\frac1z\bT-\rE\beta_1(z)\bI\bigg)^{-1}+o(\frac{1}{n})\label{mean}
\eqa where
$\beta_j=\bar{\beta}_j+\bar{\beta}_j^2\varepsilon_j+\bar{\beta}_j^2\beta_j\varepsilon_j^2$,
$\varepsilon_j=
\balpha_j'(\frac{1}{z}\bT-\bS_j)^{-1}\balpha_j-\frac1nE\rtr(\frac{1}{z}\bT-\bS_j)^{-1}$
and $\beta_j(z)=\frac{1}{1-ym(z)}+O\left(\frac1n\right)$.

By (1.15) of Bai and Silverstein (2004) and (\ref{mean}), when all
$x_{tj}$ are complex
\begin{equation}\label{mean10}R_n=0.\end{equation} In RSE case,
$$
\begin{array}{lll}
&&\frac{\rE\underline{\tilde{s}}_n(z)z}{y}\cdot R_n\nonumber\\
&=&\frac{-(\rE\beta_1(z))^2}{y}\frac{1}{n^2}\rE\rtr\bigg(\frac1z\bT-\bS_1\bigg)^{-2}{\bf T}\bigg(\frac1z\bT-\rE\beta_1(z) \bI\bigg)^{-1}+o(\frac1n)\label{R11}\\
&=&\frac{-(\rE\beta_1(z))^2}{y}\frac{1}{n^2}\rE\rtr\Bigg\{ \nonumber\\
 &&-{\bf T}\bigg(\frac1z\bT-\rE\beta_1(z) \bI\bigg)^{-1}\bigg(\rE\beta_1(z)\bI-\frac{1}{z}\bT\bigg)^{-1}\bigg(\frac1z\bT-\bS_1\bigg)^{-1}\nonumber\\
 &&-\rE\beta_1(z)\cdot\sum\limits_{k\not=1}\rE_2^{-1}(z)
                           \frac{1}{n}(\frac{1}{z}\bT-\bS_{1k})^{-1}{\bf T}\rE_1^{-1}(z)\left(\bigg(\frac1z\bT-\bS_1\bigg)^{-1}-
                           \bigg(\frac1z\bT-\bS_{1k}\bigg)^{-1}\right)\nonumber\\
 &&+\rE\beta_1(z)\cdot\sum\limits_{k\not=1}{\bf E}_2^{-1}(z)                           \bm{\alpha}_k\bm{\alpha}_k^*(\frac{1}{z}\bT-\bS_{1k})^{-1}{\bf T}{\bf E}_1^{-1}(z)\left(\bigg(\frac1z\bT-\bS_1\bigg)^{-1}-
                           \bigg(\frac1z\bT-\bS_{1k}\bigg)^{-1}\right)\nonumber\\
 &&+\rE\beta_1(z)\cdot\sum\limits_{k\not=1}{\bf E}_2^{-1}(z)
                           \left(\bm{\alpha}_k\bm{\alpha}_k^*-\frac{1}{n}{\bf I}\right)(\frac{1}{z}\bT-\bS_{1k})^{-1}
                          {\bf T} {\bf E}_1^{-1}(z)\bigg(\frac1z\bT-\bS_{1k}\bigg)^{-1}\Bigg\}+O(\frac{1}{n^{3/2}})\nonumber\\
&=&\frac{(\rE\beta_1(z))^2}{yn^2}E\rtr\bigg(\frac1z\bT-E\beta_1(z) \bI\bigg)^{-1}{\bf T}
 \bigg(E\beta_1(z)\bI-\frac{1}{z}\bT\bigg)^{-1}\bigg(\frac1z\bT-\bS_1\bigg)^{-1}\nonumber\\
&&+\frac{\rE\beta_1(z)(E\beta_1(z))^2}{yn^2}E\rtr
  \sum\limits_{k\not=1}{\bf E}_2^{-1}(z) \bm{\alpha}_k\bm{\alpha}_k^*\frac{(\frac{1}{z}\bT-\bS_{1k})^{-1}{\bf T}{\bf E}_1^{-1}(z)(\frac{1}{z}\bT-\bS_{1k})^{-1}\bm{\alpha}_k\bm{\alpha}_k^*
                     (\frac{1}{z}\bT-\bS_{1k})^{-1}}{1-\bm{\alpha}_k^*(\frac{1}{z}\bT-\bS_{1k})^{-1}\bm{\alpha}_k}
                      +O(\frac{1}{n^{3/2}})\nonumber\\
&=&\frac{(\rE\beta_1(z))^2}{yn^2}E\rtr\bigg(\frac1z\bT-E\beta_1(z) \bI\bigg)^{-1}{\bf T}
 \bigg(\rE\beta_1(z)\bI-\frac{1}{z}\bT\bigg)^{-1}\bigg(\frac1z\bT-\bS_1\bigg)^{-1}\nonumber\\
&&+y^{-1}n^{-2}(\rE\beta_1(z))^3\sum\limits_{k\not=1}E\bm{\alpha}_k^*
                     (\frac{1}{z}\bT-\bS_{1k})^{-1}{\bf T}{\bf E}_2^{-1}(z)\bm{\alpha}_k                   \cdot\frac{\bm{\alpha}_k^*(\frac{1}{z}\bT-\bS_{1k})^{-2}{\bf T}{\bf E}_1^{-1}(z)\bm{\alpha}_k}{1-\bm{\alpha}_k^*(\frac{1}{z}\bT-\bS_{1k})^{-1}\bm{\alpha}_k}
                      +O(\frac{1}{n^{3/2}})\nonumber\\
 &=&-\frac{(\rE\beta_1(z))^2}{y}\frac{1}{n^2}\rE\rtr\bigg(\frac1z\bT-E\beta_1(z) \bI\bigg)^{-1}{\bf T}
 \bigg(\rE\beta_1(z)\bI-\frac{1}{z}\bT\bigg)^{-2}\nonumber\\
 \\
&&+\frac{(\rE\beta_1(z))^4}{yn^{2}}\sum\limits_{k\not=1}\rE\bm{\alpha}_k^*
                     (\frac{1}{z}\bT-\bS_{1k})^{-1}{\bf T}{\bf E}_2^{-1}(z)\bm{\alpha}_k                   \cdot\bm{\alpha}_k^*(\frac{1}{z}\bT-\bS_{1k})^{-2}{\bf T}{\bf E}_1^{-1}(z)\bm{\alpha}_k
                     +O(\frac{1}{n^{3/2}})\nonumber\\
\end{array}
$$
where $E_1^{-1}(z)=\bigg(\frac1z\bT-E\beta_1(z) \bI\bigg)^{-1}$, $E_2^{-1}(z)=\bigg(E\beta_1(z)\bI-\frac{1}{z}\bT\bigg)^{-1}$
and
$\bA(z)=\sum\limits_{k\not=i}
                           \bigg(\frac{n-1}{n}b_i(z)\bI-\frac{1}{z}\bT\bigg)^{-1}
                           (\balpha_k\balpha_k^*-\frac{1}{n}\bI)(\frac{1}{z}\bT-\bS_{ik})^{-1}$
because
$$
\begin{array}{l}
-\frac{-1}{y(1-y\tilde{m}(z))^2}\frac{1}{n^2}E\rtr(\frac1z\bT-E\beta_1(z)
\bI)^{-1}
\sum\limits_{k\not=1}(E\beta_1(z)\bI-\frac{1}{z}\bT)^{-1}
                           \frac{1}{n}(\frac{1}{z}\bT-\bS_{1k})^{-1}\\
                           \cdot((\frac1z\bT-\bS_1)^{-1}-
                           (\frac1z\bT-\bS_{1k})^{-1})=O(\frac{1}{n^2})
\end{array}
$$
and
$$
\begin{array}{l}
\frac{-1}{y(1-y\tilde{m}(z))^2}\frac{1}{n^2}E\rtr\bigg(\frac1z\bT-E\beta_1(z)
\bI\bigg)^{-1}
  \sum\limits_{k\not=1}\bigg(E\beta_1(z)\bI-\frac{1}{z}\bT\bigg)^{-1}\\
                             \cdot\left(\bm{\alpha}_k\bm{\alpha}_k^*-\frac{1}{n}{\bf I}\right)
                           \bigg(\frac1z\bT-\bS_{1k}\bigg)^{-2}=O(\frac{1}{n^{3/2}}).
                           \end{array}
$$

%By (\ref{R1}), we have
\begin{eqnarray} && \frac{-z\rE\underline{\tilde{s}}_n(z)}{y}\cdot
pR_n=-z\rE\underline{\tilde{s}}_n(z)\cdot
nR_n\nonumber\\
&=&-\frac{\frac{(\rE\beta_1(z))^2}{y}\cdot
\frac{y}{n}\rE\rtr(\rE\beta_1(z)\bI-\frac{1}{z}\bT)^{-2}{\bf T}
\bigg(\frac1z\bT-\rE\beta_1(z) \bI\bigg)^{-1}}
{1-\frac{1}{(1-y\tilde{m}(z))^2}\cdot \frac{1}{n}\rtr
                           \bigg(\rE\beta_1(z)\bI-\frac{1}{z}\bT\bigg)^{-2}}\nonumber\\
&=&-\frac{\frac{1}{(1-y\tilde{m}(z))^2}\cdot\int
\frac{ytdH(t)}{(\rE\beta_1(z)-\frac{t}{z})^2(\frac{t}{z}-\rE\beta_1(z))}}
{1-\frac{y}{(1-y\tilde{m}(z))^2}\int\frac{dH(t)}{(E\beta_1(z)-\frac{t}{z})^2}}+o(1)\nonumber\\
&=&-\frac{\frac{1}{(1-y\tilde{m}(z))^2}\cdot\int
\frac{ytdH(t)}{(t/z-1/(1-y\tilde{m}(z)))^3}}
{1-\frac{y}{(1-y\tilde{m}(z))^2}\cdot\int\frac{dH(t)}{(t/z-1/(1-y\tilde{m}(z)))^2}}+o(1)\nonumber\\
&=&-\frac{y\int\frac{t(1+yzs(z))^2dH(t)}{(t/z-1-yzs(z))^3}}{1-y\int\frac{(1+yzs(z))^2dH(t)}{(t/z-1-yzs(z))^2}}\cdot
+o(1)\nonumber
\end{eqnarray} where $E\beta_1(z)\rightarrow\frac{1}{1-y\tilde{m}(z)}$ by
(\ref{6}). Thus by (\ref{mean10}), (\ref{tilde}) and (\ref{7}), we
have \begin{eqnarray}
p(\rE s_n(z)-s_n^0(z))&=&\frac{-1}{z^2}\cdot p(\rE\tilde{s}_n(z)-\tilde{s}_n^0(z))\nonumber\\
&=&\frac{-1}{z^2}\cdot n(\rE\underline{\tilde{s}}_n(z)-\underline{\tilde{s}}_n^0(z))\nonumber\\
&=&\frac{-1}{z^2}\cdot \frac{\rE\underline{\tilde{s}}_n(z)\underline{\tilde{s}}_n^0(z)}
{1-\frac{y\int\frac{1}
{(t+\rE\underline{\tilde{s}}_n(z))(t+\underline{\tilde{s}}_n^0(z))}dH_n(t)}
{\left(-\frac1z+y\int\frac{1}{t+\rE\underline{\tilde{s}}_n(z)}dH_n(t)-R_n\right)
\left(-\frac1z+y\int\frac{1}{t+\underline{\tilde{s}}_n^0(z)}dH_n(t)\right)}}\cdot nR_n\nonumber\\
&=&\frac{\kappa-1}{z^2}\cdot \frac{y\int
\frac{t(1+yzs(z))^3dH(t)}{(t/z-1-yzs(z))^3}}{\left(
1-y\int\frac{(1+yzs(z))^2
dH(t)}{(t/z-1-yzs(z))^2}\right)^2}+o(1)\label{bias}
\end{eqnarray}
So we conclude that in the RSE case
$$
\sup\limits_{z\in\mathcal{C}_n}\left|M_n^2(z)-\frac{\kappa-1}{z^2}\cdot
\frac{y\int \frac{t(1+yzs(z))^3dH(t)}{(t/z-1-yzs(z))^3}}{\left(
1-y\int\frac{(1+yzs(z))^2
dH(t)}{(t/z-1-yzs(z))^2}\right)^2}\right|\rightarrow0\quad\mbox{as}~n\rightarrow+\infty.
$$
By (\ref{deriv}), we obtain
$$
\rE M(z)=\frac{\kappa-1}{z^2}\cdot \frac{y\int
\frac{t(1+yzs(z))^3dH(t)}{(t/z-1-yzs(z))^3}}{\left(
1-y\int\frac{(1+yzs(z))^2 dH(t)}{(t/z-1-yzs(z))^2}\right)^2}
=\frac{-(\kappa-1)}{z^2}\cdot \frac{y\int
\frac{t(-z)^3(1+yzs(z))^3dH(t)}{(t-z-yz^2s(z))^3}}{\left(
1-y\int\frac{z^2(1+yzs(z))^2 dH(t)}{(t-z-yz^2s(z))^2}\right)^2}
$$
That is,
$$
\rE
M(z)=\frac{\kappa-1}{2}\frac{d\log\left(1-y\int\frac{(1+yzs(z))^2
dH(t)}{(t/z-1-yzs(z))^2}\right)}{dz}
$$

So the proof of Theorem \ref{thm5} is completed. \eprf

\subsection{Some Notations and Lemmas}
\begin{lem}\label{BaiLem91}(Bai and Silverstein (2010) P225) Suppose that $x_i,~i=1,\cdots,n$ are
independent, with $E x_i=0$, $E |x_i|^2=1$, $\sup E
|x_i|^4=\nu<+\infty$ and $|x_i|\leq \eta_n\sqrt{n}$ with $\eta_n>0$.
Assume that ${\bf A}$ is a complex matrix. Then, for any given
$2\leq \ell \leq b\log(n\nu^{-1}\eta_n^4)$ and $b>1$, we have
$$
E |\bm{\alpha}^*{\bf A}\bm{\alpha}-\rtr({\bf A})|^l\leq \nu
n^l(n\eta_n^4)^{-1}(40b^2\|{\bf A}\|\eta_n^2)^l
$$
where $\bm{\alpha}=(x_1,\cdots,x_n)^T$.
\end{lem}

\begin{lem}\label{2}(Bai and Silverstein (2010) P271) We have
$$
\left|E\left(\prod\limits_{k=1}^m\bgma_t^*\bA_k\bgma_t\right)
\prod\limits_{l=1}^q(\bgma_t^*\bB_l\bgma_t-n^{-1}\rtr\bT\bB_l)\right|
\leq Kn^{-(1\bigwedge
q)}\eta_n^{(2q-4)\bigvee0}\prod\limits_{k=1}^m\|\bA_k\|\prod\limits_{l=1}^q\|\bB_l\|,
$$
where $m\geq0$, $q\geq0$,
$\bgma_t=\frac{1}{\sqrt{n}}\bT^{\frac12}\bX_t$,
$\bX_t=(x_{t1},\cdots,x_{tp})^T$, $(x_{tj}, t=1,\cdots,n, j=1,\cdots,p)$ are independent with $E x_{tj}=0$, $E |x_{tj}|^2=1$,
$\sup E |x_{tj}|^4=\nu<+\infty$ and $|x_{tj}|\leq \eta_n\sqrt{n}$ with
$\eta_n>0$.
\end{lem}
\begin{lem}\label{lem3}
Under Assumptions 1-2, we obtain \be \frac1p\rtr(\frac{1}{z}{\bf
T}-{\bf S})^{-1}\rightarrow \tilde{m}(z), \ a.s. \label{eq21} \ee
where $\tilde{m}(z)$ is the unique solution to the equation
$\tilde{m}(z)=\int\frac{dH(t)}{-\frac{t}{z}\frac{1}{1-y\tilde{m}(z)}}$
satisfying
$$\Im(z)\cdot\Im(\tilde{m}(z))\ge 0.$$
\end{lem}
Proof. For any real $z<0$ and complex $w$ with $\Im(w)>0$, by
(4.1.2) of Page 61 of Bai and Silverstein (2010), we have \bqa
&&\frac1p\rtr(\frac{1}{z}{\bf T}-{\bf S}+w{\bf
I})^{-1}=-\frac1p\rtr(\frac{1}{-z}{\bf T}+{\bf S}-w{\bf
I})^{-1}\non &\rightarrow& -\tilde
m(z,w)=-m_{H_z}\left(w-\frac{1}{y}\int\frac{\tau
dH_0(\tau)}{1+\tau \tilde m(z,w)}\right),\ a.s. \label{eq22} \eqa
where $\tilde{m}(z,w)$ is limit of the Stieltjes transform of the
matrix $\frac{1}{-z}{\bf T}+{\bf S}$, $m_{H_z}$ is the Stieltjes
transform of $H_z$, the LSD of $\frac1{-z}{\bf T}$, and
$H_0(\tau)=I_{(\tau>y)}$. By Theorem 5.11 and Lemma 2.14 (Vitali
Lemma) of Bai and Silverstein (2010), the convergence of
(\ref{eq22}) is also true for $w=0$.
 That is,
\bqa &&\frac1p\rtr(\frac{1}{z}{\bf T}-{\bf S})^{-1}\to -\tilde
m(z,0)=-m_{H_z}\left(-\frac{1}{1+y\tilde m(z,0)}\right),\ a.s.\non
&=&-\int\frac{1}{\lambda-\frac{1}{1-y\tilde
m(z,0)}}dH_z(\lambda)=-
\int\frac{1}{\frac{\lambda}{-z}+\frac{1}{1+y\tilde
m(z,0)}}dH(\lambda),\ a.s. \label{eq23} \eqa Denoting
$\tilde{m}(z)=-\tilde m(z,0)$, then the convergence of
(\ref{eq21}) is proved for all real nonpositive $z$. Noting that
both sides of (\ref{eq21}) are analytic functions of $z$ on the
region $D^-=\{z\in \mathbb C: z\mbox{ is not nonpositive real
number}\}$, applying Vitali Lemma again, we conclude that
(\ref{eq21}) is true for all $z\in D^-$ and $\tilde{m}(z)$
satisfies
\begin{equation}\label{eq24}
\tilde{m}(z)=\int\frac{1}{\frac{\lambda}{z}-\frac{1}{1-y
\tilde{m}(z)}}dH(\lambda)
\end{equation}
Because the imaginary part of LHS of (\ref{eq24}) has the same
sign as $z$, we conclude that $\Im (\tilde{m}(z))$ should have the
same sign as $\Im(z)$.

Our next goal is to show that for every non-real $z$, the equation
(\ref{eq24}) has a unique solution $\tilde{m}(z)$ whose imaginary
part has the same sign as $\Im(z)$. By symmetry, we only need to
consider the case where $\Im(z)>0$. Suppose that there are two
different solutions $m_1(z)\ne m_2(z)$. Making difference of both
sides of (\ref{eq24}), we obtain \bqa 1&=&\int \frac{\frac{y}{(1-y
m_1)(1-ym_2)}}{(\frac{\lambda}{z}-\frac{1}{1-y
m_1})(\frac{\lambda}{z}-\frac{1}{1-y m_2})}dH(\gl)\non
&\le&\left(\int \frac{\frac{y}{|1-y
m_1|^2}}{\left|\frac{\lambda}{z}-\frac{1}{1-y
m_1}\right|^2}dH(\gl)\int \frac{\frac{y}{|1-y
m_2|^2}}{\left|\frac{\lambda}{z}-\frac{1}{1-y
m_2}\right|^2}dH(\gl)\right)^{1/2}. \label{eq25} \eqa Comparing
the imaginary parts of both sides of (\ref{eq25}), we have \bqn
\Im(m_j)=\int
\frac{\frac{\Im(z)\lambda}{|z|^2}+\frac{y\Im(m_j)}{|1-y
m_j|^2}}{\left|\frac{\lambda}{z}-\frac{1}{1-y
m_j}\right|^2}dH(\gl),\ j=1,2. \eqn Since $\Im(m_j)>0$ implies
that
$$
\int \frac{\frac{y}{|1-y m_j|^2}}{\left|\frac{\lambda}{z}-\frac{1}{1-y m_j}\right|^2}dH(\gl)<1,
$$
which contradicts to (\ref{eq25}).

The proof of the lemma is completed. \eprf

\begin{lem}\label{lem4}
Under Assumptions 1, 2, 3, we have
$$
s(z)=-\frac1z-\frac{1}{z^2}\cdot\tilde{s}(z)
$$
where
$$
\underline{\tilde{s}}(z)=\frac{-z}{1-y\tilde{m}(z)}=\frac{-z}{1-y\int\frac{1}{t/z-(1-y\tilde{m}(z))^{-1}}dH(t)}
=\frac{1}{-\frac1z+y\int\frac{1}{t+\underline{\tilde{s}}(z)}dH(t)},
$$
$\underline{\tilde{s}}(z)=-(1-y)z+y\tilde{s}(z)$, $s(z)$ is the
Stieltjes transform of the LSD of ${\bf S}^{-1}{\bf T}$, and
$\tilde{s}(z)$ is the limit of $\frac1p\rtr\left({\bf
S}-\frac{{\bf T}}{z}\right)^{-1}{\bf T}$.

\proof We have
$$
\frac1p\rtr({\bf S}^{-1}{\bf T}-z{\bf I})^{-1}=-\frac1z-\frac{1}{z^2}\cdot\frac1p\rtr\left({\bf S}-\frac{{\bf T}}{z}\right)^{-1}{\bf T}
$$
and
$$
s_n(z)=\frac{1}{pz}\sum\limits_{i=1}^n\frac{\balpha_i'(\frac{1}{z}\bT-\bS_i)^{-1}\balpha_i}
         {1-\balpha_i'(\frac{1}{z}\bT-\bS_i)^{-1}\balpha_i}
=-\frac{1}{yz}+\frac{1}{yz}\cdot\frac1n\sum\limits_{i=1}^n\frac{1}
         {1-\balpha_i'(\frac{1}{z}\bT-\bS_i)^{-1}\balpha_i}
$$
where $\balpha_i=\frac{1}{\sqrt{n}}\bX_i$, $i=1\cdots,n$.
Let the limit of $\frac1p\rtr\left({\bf S}-\frac{{\bf T}}{z}\right)^{-1}{\bf T}$
be $\tilde{s}(z)$ and $\tilde{s}_n(z)=\frac1p\rtr\left({\bf S}-\frac{{\bf T}}{z}\right)^{-1}{\bf T}$. Let
$$
\underline{\tilde{s}}(z)=-(1-y)z+y\tilde{s}(z).
$$
In fact, we have
\begin{equation}\label{tilde}
s_n(z)=-\frac1z-\frac{1}{z^2}\tilde{s}_n(z)
\end{equation}
$$
s(z)=-\frac1z-\frac{1}{z^2}\cdot\tilde{s}(z)
=-\frac{1}{yz}+\frac{1}{yz}\cdot\frac{1}{1-y\tilde{m}(z)}
$$
\begin{equation}\label{tilde1}
\underline{\tilde{s}}(z)=-(1-y)z+y\tilde{s}(z)=\frac{-z}{1-y\tilde{m}(z)},\quad
E\underline{\tilde{s}}_n(z)=-zE\beta_i(z)
\end{equation}
where
$\beta_i(z)=\frac{1}{1-\bm{\alpha}_i'(\frac{\bT}{z}-\bS_i)^{-1}\bm{\alpha}_i}$.
Therefore, we have
$$
\underline{\tilde{s}}(z)=\frac{-z}{1-y\tilde{m}(z)}=\frac{-z}{1-y\int\frac{1}{\frac{t}{z}-\frac{1}{1-y\tilde{m}(z)}}dH(t)}
=\frac{1}{-\frac1z+y\int\frac{1}{t+\underline{\tilde{s}}(z)}dH(t)}.
$$
 That is,
$$
\underline{\tilde{s}}(z)=\frac{1}{-\frac1z+y\int\frac{1}{t+\underline{\tilde{s}}(z)}dH(t)}.
$$

The proof of the lemma is completed. \eprf
\end{lem}
Here we give some notes:
\begin{equation}
\tilde{m}(z)=\int\frac{1}{\frac{\lambda}{z}-\frac{1}{1-y
\tilde{m}(z)}}dH(\lambda)
\end{equation}
where $\tilde{m}(z)$ is the limit of
$\frac1p\rtr(\frac1z\bbT-\bbS)^{-1}$ and $H(t)$ is the LSD of
$\bbT$.
\begin{equation}\label{7}
\underline{\tilde{s}}_n^0(z)\rightarrow\underline{\tilde{s}}(z)=\frac{-z}{1-y\tilde{m}(z)},\quad
1+yzs(z)=\frac{1}{1-y\tilde{m}(z)}=\frac{1}{1-y\int\frac{dH(t)}{t/z-(1+yzs(z))}}
\end{equation}
where $s(z)$ is the Stieltjes transform of the LSD of
$\bbS^{-1}\bbT$ and
$\underline{\tilde{s}}_n^0(z)=\frac{1}{-\frac1z+y\int\frac{1}{t+\underline{\tilde{s}}_n^0(z)}dH_n(t)}$
with the ESD $H_n(t)$ of $\bbT$.
\begin{equation}\label{6}
E\beta_1(z)\rightarrow \frac{1}{1-y\tilde{m}(z)}=1+yzs(z),\quad
b_i(z)\rightarrow \frac{1}{1-y\tilde{m}(z)}=1+yzs(z)
\end{equation}
where
$\beta_1(z)=\frac{1}{1-\bm{\alpha}_1'(\frac{\bT}{z}-\bS_1)^{-1}\bm{\alpha}_1}$
and $b_1(z)=\frac{1}{1-n^{-1}\rE\rtr(\frac{1}{z}\bT-\bS_1)^{-1}}$.
\begin{equation}\label{8}
-z(1+yzs(z))=\frac{-z}{1-y\int\frac{dH(t)}{t/z-(1+yzs(z))}}=\frac{-z}{1-y\int\frac{zdH(t)}{t-z(1+yzs(z))}}.
\end{equation}
\begin{equation}\label{9}
-\frac{1}{z}+y\int\frac{dH(t)}{t-z(1+yzs(z))}=\frac{1}{-z(1+yzs(z)}
\end{equation}
\begin{equation}\label{10}
\frac{1}{z^2}-y\int\frac{(-z(1+yzs(z)))'dH(t)}{(t-z(1+yzs(z)))^2}=\frac{-(-z(1+yzs(z))'}{(-z(1+yzs(z))^2}
\end{equation}
\begin{equation}\label{11}
\frac{(z(1+yzs(z))^2}{z^2}-y\int\frac{(z(1+yzs(z)))^2dH(t)}{(t-z(1+yzs(z)))^2}(-z(1+yzs(z))'
=-(-z(1+yzs(z))'
\end{equation}
\begin{equation}\label{deriv}
(-z(1+yzs(z))'=\frac{-1}{z^2}\cdot\frac{(-z(1+yzs(z))^2}{1-y\int\frac{(-z(1+yzs(z)))^2dH(t)}{(t-z(1+yzs(z)))^2}}.
\end{equation}

Especially, when $\bbT=\bbI_p$, by (\ref{sz}), (\ref{7}),
(\ref{eq24}) and the definition of $\tilde{m}(z)$, we have
$$
1+yzs(z)=\frac{1}{1+y\cdot
m(\frac1z)}=-\frac1z\underline{m}(\frac1z)=(1-y)-\frac{y}{z}\cdot
m(\frac1z)
$$
and
$$\tilde{m}(z)=-m(\frac1z)=\frac{1}{\frac{1}{z}-\frac{1}{1+y\cdot
m(\frac1z))}}.$$

\end{document}